\documentclass[aos,MSNbibl,seceqn,nameyear,dvips]{arximspdf}
\usepackage{dcolumn}

\doi{10.1214/13-AOS1159} 
\volume{41}
\issue{5}
\pubyear{2013}
\firstpage{2505}
\lastpage{2536}

\makeatletter
\newcolumntype{d}[1]{D{.}{.}{#1}}
\newcommand{\rrVert}{\Vert}
\newcommand{\llVert}{\Vert}
\newcommand{\rrvert}{\vert}
\newcommand{\llvert}{\vert}
\newproclaim{rem}{Remark}
\newproclaim{definition}{Definition}[section]
\newtheorem{lemma}[definition]{Lemma}
\newtheorem{theorem}[definition]{Theorem}
\newtheorem{corollary}[definition]{Corollary}
\newproclaim{exa}{Example}
\makeatother

\begin{document}
\begin{frontmatter}

\title{Calibrating nonconvex penalized regression in~ultra-high dimension}
\runtitle{Nonconvex penalized regression}

\begin{aug}
\author[A]{\fnms{Lan} \snm{Wang}\thanksref{t1}\ead[label=e1]{wangx346@umn.edu}},
\author[B]{\fnms{Yongdai} \snm{Kim}\thanksref{t2}\ead[label=e2]{ydkim0903@gmail.com}}
\and
\author[C]{\fnms{Runze} \snm{Li}\corref{}\thanksref{t3}\ead[label=e3]{rzli@psu.edu}}
\runauthor{L. Wang, Y. Kim and R. Li}
\affiliation{University of Minnesota, Seoul National University and\\
Pennsylvania State University}
\address[A]{L. Wang\\
School of Statistics\\
University of Minnesota\\
Minneapolis, Minnesota 55455\\
USA\\
\printead{e1}}
\address[B]{Y. Kim\\
Department of Statistics \\
Seoul National University\\
Seoul, Korea\\
\printead{e2}}
\address[C]{R. Li\\
Department of Statistics\\
\quad and the Methodology Center\\
Pennsylvania State University\\
University Park, Pennsylvania 16802\\
USA\\
\printead{e3}}

\end{aug}
\thankstext{t1}{Supported in part by NSF Grant DMS-13-08960.}
\thankstext{t2}{Supported in part by National Research Foundation of
Korea Grant number 20100012671 funded by the Korea government.}
\thankstext{t3}{Supported in part by National Natural Science
Foundation of China, 11028103 and NIH Grants
P50 DA10075, R21 DA024260, R01 CA168676 and R01 MH096711.}

\received{\smonth{9} \syear{2012}}
\revised{\smonth{6} \syear{2013}}

%
\begin{abstract}
We investigate high-dimensional nonconvex penalized regression, where
the number of covariates may grow at an exponential rate. Although
recent asymptotic theory established that there exists a local minimum
possessing the oracle property under general conditions, it is still
largely an open problem how to identify the oracle estimator among
potentially multiple local minima. There are two main obstacles: (1)
due to the presence of multiple minima, the solution path is nonunique
and is not guaranteed to contain the oracle estimator; (2) even if a
solution path is known to contain the oracle estimator, the optimal
tuning parameter depends on many unknown factors and is hard to
estimate. To address these two challenging issues, we first prove that
an easy-to-calculate calibrated CCCP algorithm produces a consistent
solution path which contains the oracle estimator with probability
approaching one. Furthermore, we propose a high-dimensional BIC
criterion and show that it can be applied to the solution path to
select the optimal tuning parameter which asymptotically identifies the
oracle estimator. The theory for a general class of nonconvex penalties
in the ultra-high dimensional setup is established when the random
errors follow the sub-Gaussian distribution. Monte Carlo studies
confirm that the calibrated CCCP algorithm combined with the proposed
high-dimensional BIC has desirable performance in identifying the
underlying sparsity pattern for high-dimensional data analysis.
\end{abstract}

%

\end{frontmatter}
\newpage

\section{Introduction}\label{sec1}

High-dimensional data, where the number of covariates $p$~greatly
exceeds the sample size $n$, arise frequently in modern applications in
biology, chemometrics, economics, neuroscience and other scientific
fields. To facilitate the analysis, it is often useful and reasonable
to assume that only a small number of covariates are relevant for
modeling the response variable. Under this sparsity assumption, a
widely used approach for analyzing high-dimensional data is regularized
or penalized regression. This approach estimates the unknown regression
coefficients by solving the following penalized regression problem:
%
%
%
%
\begin{equation}
\label{model2} \min_{\bolds{\beta}\in\mathcal
{R}^p} \Biggl\{ (2n)^{-1}\|
\mathbf{y}-\mathbf{X} \bolds{\beta}\|^2+\sum
_{j=1}^pp_{\lambda}\bigl(\llvert
\beta_j\rrvert \bigr) \Biggr\},
\end{equation}
where $\mathbf{y}$ is the vector of responses, $\mathbf{X}$ is an
$n\times p$ matrix of covariates,
$\bolds{\beta}=(\beta_1,\ldots,\beta_p)^T$ is the vector of unknown
regression coefficients, $\|\cdot\|$ denotes the $L_2$ norm (Euclidean
norm), and $p_{\lambda}(\cdot)$ is a penalty function which depends on
a tuning parameter $\lambda>0$. Many commonly used variable selection
procedures in the literature can be cast into the above framework,
including the best subset selection, $L_1$ penalized regression or
Lasso [\citet{T1}], Bridge regression [\citet{Fr}], SCAD
[\citet{F1}], MCP [\citet{Zh1}], among others.

The Lasso penalized regression is computationally attractive and enjoys
great performance in prediction. However, it is known that Lasso
requires rather stringent conditions on the design matrix to be
variable selection consistent [\citet{Z1}, \citet{Zhao1}].
Focusing on identifying the unknown sparsity pattern, nonconvex
penalized high-dimensional regression has recently received
considerable attention. \citet{F1} first systematically studied
nonconvex penalized likelihood for fixed finite dimension $p$. In
particular, they recommended the SCAD penalty which enjoys the oracle
property for variable selection. That is, it can estimate the zero
coefficients as exact zero with probability approaching one, and
estimate the nonzero coefficients as efficiently as if the true
sparsity pattern is known in advance. \citet{F3} extended these
results by allowing $p$ to grow with $n$ at the rate $p=o(n^{1/5})$ or
$p=o(n^{1/3})$. For high dimensional nonconvex penalized regression
with $p\gg n$, \citet{K1} proved that the oracle estimator itself
is a local minimum of SCAD penalized least squares regression under
very relaxed conditions; \citet{Zh1} proposed a minimax concave
penalty (MCP) and devised a novel PLUS algorithm which when used
together can achieve the oracle property under certain regularity
conditions. Important insight has also been gained through the recent
work on theoretical analysis of the global solution
[\citet{K20,Zh3}]. However, direct computation of the global
solution to the nonconvex penalized regression is infeasible in high
dimensional setting.

For practical data analysis, it is critical to find an
easy-to-implement procedure which can find a local solution with
satisfactory theoretical property even when the number of covariates
greatly exceeds the sample size. Two challenging issues remain
unsolved. One is the problem of multiple local minima; the other is
the problem of optimal tuning parameter selection.

A direct consequence of the multiple local minima problem is that the
solution path is not unique and is not guaranteed to contain the oracle
estimator. This problem is due to the nature of the nonconvexity of the
penalty. To understand it, we note that the penalized objective
function in (\ref{model2}) is nonconvex in $\bolds{\beta}$ whenever the
convexity of the least squares loss function does not dominate the
concavity of the penalty part. In general, the occurrence of multiple
minima is unavoidable unless strong assumptions are imposed on both the
design matrix and the penalty function. The recent theory for SCAD
penalized linear regression [\citet{K1}] and for general
nonconcave penalized generalized linear models [\citet{F2}]
indicates that one of the local minima enjoys the oracle property but
it is still an unsolved problem how to identify the oracle estimator
among multiple minima when $p\gg n$. Popularly used algorithms
generally only ensure the convergence to a local minimum, which is not
necessarily the oracle estimator. Numerical evidence in
Section~\ref{sec4} suggests that the local minima identified by some of
the popular algorithms have a relatively low probability to recover the
unknown sparsity pattern although it may have small estimation error.

Even if a solution path is known to contain the oracle estimator,
identifying such a desirable estimator from the path is itself a
challenging problem in ultra-high dimension. The main issue is to find
the optimal tuning parameter which yields the oracle estimator. The
theoretically optimal tuning parameter does not have an explicit
representation and depends on unknown factors such as the variance of
the unobserved random noise. Cross-validation is commonly adopted in
practice to select the tuning parameter but is observed to often result
in overfitting. In the case of fixed $p$, \citet{W2} rigorously
proved that generalized cross-validation leads to an overfitted model
with a positive probability for SCAD-penalized regression. Effective
BIC-type criterion for nonconvex penalized regression has been
investigated in \citet{W2} and \citet{Zh4} for fixed $p$; and
in \citet{W1} for diverging $p$ (but $p<n$). However, to the best
of our knowledge, there is still no satisfactory tuning parameter
selection procedure for nonconvex penalized regression in ultra-high
dimension.

The above two main concerns motivate us to consider calibrating
nonconvex penalized regression in ultra-high dimension with the goal to
identify the oracle estimator with high probability. To achieve this,
we first prove that a calibration of the CCCP algorithm
[\citet{K1}] for nonconvex penalized regression produces a
consistent solution path with probability approaching one in merely two
steps under conditions much more relaxed than what would be required
for the Lasso estimator to be model selection consistent. Furthermore,
extending the recent work of \citet{C2} and \citet{K2} for
Bayesian information criterion (BIC) on high dimensional least squares
regression, we propose a high-dimensional BIC for a nonconvex penalized
solution path and prove its validity under more general conditions when
$p$ grows at an exponential rate. The recent independent work of
\citeauthor{Zh1} (\citeyear{Tong1,Tong2}) devised a multi-stage convex
relaxation scheme and proved that for the capped $L_1$ penalty the
algorithm can find a consistent solution path with probability
approaching one under certain conditions. Despite the similar flavor
shared with the algorithm proposed in this paper, his algorithm takes
multiple steps (which can be very large in practice depending on the
design condition) and the paper has not studied the problem of tuning
parameter selection.

To deepen our understanding of the nonconvex penalized regression, we
also derive an interesting auxiliary theoretical result of an upper
bound on the $L_2$ distance between a sparse local solution of
nonconvex penalized regression and the oracle estimator. This result is
new and insightful. It suggests that under general regularity
conditions a sparse local minimum can often have small estimation error
even though it may not be the oracle estimator. Overall, the
theoretical results in this paper fill in important gaps in the
literature, thus substantially enlarge the scope of applications of
nonconvex penalized regression in ultra-high dimension. In Monte Carlo
studies, we demonstrate that the calibrated CCCP algorithm combined
with the proposed high-dimensional BIC is effective in identifying the
underlying sparsity pattern.

The rest of the paper is organized as follows. In Section~\ref{sec2},
we define the notation, review the CCCP algorithm and introduce the new
methodology. In Section~\ref{sec3}, we establish that the proposed
calibrated CCCP solution path contains the oracle estimator with
probability approaching one under general conditions, and that the
proposed high-dimensional BIC is able to select the optimal tuning
parameter with probability tending to one. In Section~\ref{sec4}, we
report numerical results from Monte Carlo simulations and a real data
example. In Section~\ref{sec5}, we present an auxiliary theoretical
result which sheds light on the estimation accuracy of a local minimum
of nonconvex penalized regression if it is not the oracle estimator.
The proofs are given in Section~\ref{sec6}.

\section{Calibrated nonconvex penalized least squares method}\label{sec2}

\subsection{Notation and setup}\label{sec2.1}

Suppose that $\{(Y_i,\mathbf{x}_i)\}_{i=1}^n$ is a random sample from
the linear regression model
%
%
%
%
\begin{equation}
\label{model1} \mathbf{y}=\mathbf{X}\bolds{\beta}^*+\bolds{\varepsilon},
\end{equation}
where $\mathbf{y}=(Y_1,\ldots,Y_n)^T$, $\mathbf{X}$ is the $n\times p$
nonstochastic design matrix with the $i$th row $\mathbf{x}_i^T$,
$\bolds{\beta}^*=(\beta_1^*,\ldots,\beta_p^*)^T$ is the vector of
unknown true parameters, and
$\bolds{\varepsilon}=(\varepsilon_1,\ldots,\varepsilon_n)^T$ is a
vector of independent and identically distributed random errors.

We are interested in the case where $p=p_n$ greatly exceeds the sample
size $n$. The vector of the true parameters $\bolds{\beta}^*$ is
assumed to be sparse in the sense that the majority of its components
are exactly zero. Let $ A_0=\{j\dvtx \beta_j^*\neq0\}$ be the index set
of covariates with nonzero coefficients and let $|A_0|=q$ denote the
cardinality of $A_0$. We use $d_{*}=\min\{|\beta_j^*|\dvtx
\beta_j^*\neq0\}$\label{page5} to denote the minimal absolute value of
the nonzero coefficients. Without loss of generality, we may assume
that the first $q$ components of $\bolds{\beta}^*$ are nonzero, thus we
can write $\bolds{\beta}^*=(\bolds{\beta}_1^{*T},\mathbf{0}^T)^T$,
where $\mathbf{0} $ represents a zero vector of length $p-q$. The
oracle estimator is defined as
$\widehat{\bolds{\beta}}{}^{(o)}=(\widehat{\bolds{\beta
}}{}_1^{(o)T},\mathbf{0}^T)^T$, where
$\widehat{\bolds{\beta}}{}_1^{(o)}$ is the least squares estimator
fitted using only the covariates whose indices are in $A_0$.

To handle the high-dimensional covariates, we consider the penalized
regression in (\ref{model2}). The penalty function $p_{\lambda}(t)$ is
assumed to be increasing and concave for $t\in[0,+\infty)$ with a
continuous derivative $\dot{p}_{\lambda}(t)$ on $(0,+\infty)$. To
induce sparsity of the penalized estimator, it is generally necessary
for the penalty function to have a singularity at the origin, that is,
$\dot{p}_{\lambda}(0+)>0$. Without loss of generality, the penalty
function can be standardized such that $\dot{p}_{\lambda}(0+)=\lambda$.
Furthermore, it is required that
%
%
%
%
\begin{eqnarray}
\label{pc1} \dot{p}_{\lambda}(t)&\leq& \lambda\qquad\forall0<t<a_0\lambda,
\\
\label{pc2} \dot{p}_{\lambda}(t)&=&0\qquad\forall t>a_0\lambda
\end{eqnarray}
for some positive constant $a_0$. Condition (\ref{pc2}) plays the key
role of not over-penalizing large coefficients, thus alleviating the
bias problem associated with Lasso.

The above class of penalty functions include the popularly used SCAD penalty
and MCP. The SCAD penalty is defined by
%
%
%
%
\begin{equation}
\label{SCAD} \dot{p}_{\lambda}(t)=\lambda \biggl\{I(t\leq\lambda)+
\frac{(a\lambda-t)_{+}}{(a-1)\lambda}I(t>\lambda) \biggr\}
\end{equation}
for some $a>2$, where the notation $b_{+}$ stands for the positive part
of $b$, that is, $b_{+}=bI(b>0)$. \citet{F1} recommended to use
$a=3.7$ from a Bayesian perspective. On the other hand, the MCP is
defined by $ \dot{p}_{\lambda}(t)=a^{-1} (a\lambda-t)_{+} $ for some
$a>0$ (as $a\downarrow1$, it amounts to hard-thresholding, thus in the
following we assume $a>1$).

Let $\mathbf{x}_{(j)}$ be the $j$th column vector of $\mathbf{X}$.
Without loss of generality, we assume that
$\mathbf{x}_{(j)}^{T}\mathbf{x}_{(j)}/n=1$ for all $j$. Throughout\vspace*{1pt} this
paper, the following notation is used. For an arbitrary index set
$A\subseteq\{1,2,\ldots,p\}$, $\mathbf{X}_A$ denotes the $n\times|A|$
submatrix of $\mathbf{X}$ formed by those columns of $\mathbf{X}$ whose
indices are in $A$. For a~vector $\mathbf{v}=(v_1,\ldots,v_p)'$, we use
$\|\mathbf{v}\|$ to denote its $L_2$ norm; on the other hand
$\|\mathbf{v}\|_0=\#\{j\dvtx v_j\neq0\}$ denotes the $L_0$ norm,
$\|\mathbf{v}\|_1=\sum_{j}|v_j|$ denotes the $L_1$ norm and
$\|\mathbf{v} \|_{\infty}=\max_j |v_j|$ denotes the $L_\infty$ norm. We
use $\mathbf{v} _{A}$ to represent the size-$|A|$ subvector of
$\mathbf{v}$ formed by the entries $v_j$ with indices in $A$. For a
symmetric matrix $\mathbf{B}$, $\lambda_{\min}(\mathbf{B})$ and
$\lambda_{\max}(\mathbf{B})$ stand for the smallest and largest
eigenvalues of $\mathbf{B}$, respectively. Furthermore, we let
%
%
%
%
\begin{equation}
\label{ximin} \xi_{\min}(m)=\min_{|B|\leq m,
A_0\subseteq B}
\lambda_{\min} \bigl(n^{-1}\mathbf{X}_B^T
\mathbf{X}_B \bigr).
\end{equation}
Finally, $p$, $q$, $\lambda$ and other related quantities are all
allowed to depend on $n$, but we suppress such dependence for
notational simplicity.

\subsection{The CCCP algorithm}\label{sec2.2}
It is challenging to solve the penalized regression problem in
(\ref{model2}) when the penalty function is nonconvex. \citet{K1}
proposed a fast optimization algorithm called the SCAD--CCCP (CCCP
stands for ConCave Convex procedure) algorithm for solving the
SCAD-penalized regression. The key idea is to update the solution with
the minimizer of the tight convex upper bound of the objective function
obtained at the current solution. What makes a fast algorithm practical
relies on the possibility of decomposing the nonconvexed penalized
least squares objective function as the sum of a convex function and a
concave function. To be specific, suppose we want to minimize an
objective function $C(\bolds{\beta})$ which has the representation
$C(\bolds{\beta})=C_{\mathrm{vex}}(\bolds{\beta}
)+C_{\mathrm{cav}}(\bolds{\beta})$ for a convex function
$C_{\mathrm{vex}}(\bolds{\beta})$ and a concave function
$C_{\mathrm{cav}}(\bolds{\beta})$. Given a current solution
$\bolds{\beta}^{(k)}$, the tight convex upper bound of
$C(\bolds{\beta})$ is given by
$Q(\bolds{\beta})=C_{\mathrm{vex}}(\bolds{\beta})+\nabla
C_{\mathrm{cav}}(\bolds {\beta}^{(k)})^{\prime}\bolds{\beta}$ where
$\nabla C_{\mathrm{cav}}(\bolds{\beta})=\partial
C_{\mathrm{cav}}(\bolds{\beta })/\partial\bolds{\beta}$. We then update
the solution by minimizing $Q(\bolds{\beta})$. Since $Q(\bolds{\beta})$
is a convex function, it can be easily minimized.

For the penalized regression in (\ref{model2}), we consider a penalty
function $p_{\lambda}(\llvert \beta_j\rrvert )$ which has the
decomposition
%
%
%
%
\begin{equation}
\label{decom} p_\lambda\bigl(\llvert \beta_j\rrvert
\bigr)=J_\lambda \bigl(| \beta_j| \bigr)+\lambda| \beta_j|,
\end{equation}
where $J_\lambda(\llvert \beta_j\rrvert )$ is a differentiable concave
function. For example, for the SCAD penalty,
\begin{eqnarray*}
J_{\lambda}\bigl(\llvert \beta_j\rrvert \bigr)
&=& -\frac{\beta_j^2-2\lambda|\beta_j|+\lambda^2}{2(a-1)} I \bigl(\lambda\leq |\beta_j| \leq a\lambda\bigr)
\\
&&{}  + \biggl[\frac{(a+1)\lambda^2}{2}-\lambda|\beta_j| \biggr] I \bigl(\llvert
\beta_j| > a\lambda\bigr),
\end{eqnarray*}
while for the MCP penalty,
\[
J_{\lambda}\bigl(\llvert \beta_j\rrvert \bigr)=
\frac{\beta_j^2}{2a}I \bigl(0\leq| \beta_j| < a\lambda\bigr) + \biggl[
\frac{a\lambda^2}{2}-\lambda| \beta_j| \biggr]I \bigl(\llvert
\beta_j| \geq a\lambda\bigr).
\]
Hence, using the decomposition in (\ref{decom}),
the penalized objective function in (\ref{model2}) can be rewritten as
\[
\frac{1}{2n}\|\mathbf{y}-\mathbf{X}\bolds{\beta}
\|^2+\sum_{j=1}^p
J_\lambda\bigl(\llvert \beta_j\rrvert \bigr)+ \lambda\sum
_{j=1}^p|\beta_j|,
\]
which is the
sum of convex and concave functions. The CCCP algorithm is applied as
follows. Given a current solution $\bolds{\beta}^{(k)}$, the tight convex
upper bound is
%
%
%
%
\begin{equation}
\label{Qfun} Q \bigl(\bolds{\beta}\mid \bolds{\beta}^{(k)}, \lambda \bigr)=
\frac{1}{2n}\|\mathbf{y}-\mathbf{X}\bolds{\beta}\|^2+ \sum
_{j=1}^p \nabla J_\lambda \bigl(
\bigl| \beta_j^{(k)}\bigr| \bigr) \beta_j + \lambda
\sum_{j=1}^p|\beta_j|.
\end{equation}
We then update the current solution by $\bolds{\beta}^{(k+1)}=\arg
\min_{\bolds{\beta}}Q(\bolds{\beta}\mid \bolds{\beta}^{(k)}, \lambda)$.

An important property of the CCCP algorithm is that the objective
function always decreases after each iteration [\citet{YR}, and
\citet{T2}], from which it can be deduced that the solution
converges to a local minimum. See, for example, Corollary~3.2 of
\citet{H3}. However, there is no guarantee that the local minimum
found is the oracle estimator itself because there are multiple local
minima and the solution of the CCCP algorithm depends on the choice of
the initial solution.

\subsection{Calibrated nonconvex penalized regression}\label{sec2.3}

In this paper, we propose and study a calibrated CCCP estimator. More
specifically, we start with the initial value
$\bolds{\beta}^{(0)}=\mathbf{0}$ and a tuning parameter $\lambda>0$
and let $Q$ be the tight convex upper bound defined in (\ref{Qfun}).
The calibrated algorithm consists of the following two steps.
\begin{enumerate}[1.]
\item Let $\widehat{\bolds{\beta}}{}^{(1)}(\lambda)=\arg
    \min_{\bolds{\beta}}Q(\bolds{\beta}\mid \bolds{\beta}^{(0)},
    \tau\lambda)$, where the choice of $\tau>0$ will be discussed
    later.
\item Let $\widehat{\bolds{\beta}}(\lambda)=\arg\min_{\bolds
    {\beta}}Q(\bolds{\beta}
    \mid\widehat{\bolds{\beta}}{}^{(1)}(\lambda), \lambda)$.
\end{enumerate}
When we consider a sequence of tuning parameter values, we obtain a
solution path $\{\widehat{\bolds{\beta}}(\lambda)\dvtx \lambda>0\}$.
The calculation of the path is fast even for very high-dimensional $p$
as for each of the two steps a convex minimization problem is solved.
In step~1, a smaller tuning parameter $\tau\lambda$ is adopted to
increase the estimation accuracy, see Section~\ref{sec3.1} for
discussions on the practical choice of $\tau$. We call a solution path
``\textit{path consistent}'' if it contains the oracle estimator. In
Section~\ref{sec3.1}, we will prove that the calibrated CCCP algorithm
produces a consistent solution path under rather weak conditions.

Given such a solution path, a critical question is how to tune the
regularization parameter $\lambda$ in order to identify the oracle
estimator. The performance of a penalized regression estimator is known
to heavily depend on the choice of the tuning parameter. To further
calibrate nonconvex penalized regression, we consider the following
high-dimensional BIC criterion (HBIC) to compare the estimators from
the above solution path:
%
%
%
%
\begin{equation}
\label{hbic} \mathrm{HBIC}(\lambda)=\log \bigl(\widehat{\sigma}{}^2_{\lambda
} \bigr)+|M_{\lambda}|\frac{C_n\log(p)}{n},
\end{equation}
where $M_{\lambda}=\{j\dvtx \widehat{\bolds{\beta}}_j(\lambda)\neq0\} $
is the model identified by $\widehat{\bolds{\beta}}(\lambda)$,
$|M_{\lambda}|$ denotes the cardinality of $M_{\lambda}$, and
$\widehat{\sigma}{}^2_{\lambda}=n^{-1}\mathrm{SSE}_{\lambda}$ with
$\mathrm{SSE}_{\lambda}=\|\mathbf{Y}-\mathbf{X}\widehat{\bolds{\beta}}(\lambda)\|^2$.
As we are interested in the case where $p$ greatly exceeds $n$, the
penalty term also depends on~$p$; and $C_n$ is a sequence of numbers
that diverges to $\infty$, which will be discussed later.

We compare the value of the above HBIC criterion for $\lambda\in
\Lambda_n=\{\lambda\dvtx |M_{\lambda}|\leq K_n\}$, where $K_n>q$
represents a rough estimate of an upper bound of the sparsity of the
model and is allowed to diverge to $\infty$. We select the tuning
parameter
\[
\widehat{\lambda}=\mathop{\arg\min}_{\lambda\in\Lambda_n}\mathrm{HBIC}(\lambda).
\]

The above criterion extends the recent works of \citet{C2} and
\citet{K2} on the high-dimensional BIC for the least squares
regression to tuning parameter selection for nonconvex penalized
regression. In Sections~\ref{sec3.1}--\ref{sec3.3}, we study asymptotic
properties under conditions such as sub-Gaussian random errors,
dimension of the covariates growing at the exponential rate and
diverging~$K_n$.

\section{Theoretical properties}\label{sec3}

The main theory comprises two parts. We first show that under some
general regularity conditions the calibrated CCCP algorithm yields a
solution path with the ``\textit{path consistency}'' property. We next
verify that when the proposed high-dimensional BIC is applied to this
solution path to choose the tuning parameter $\lambda$, with
probability tending to one the resulted estimator is the oracle
estimator itself.

To facilitate the presentation, we specify a set of regularity conditions.
\begin{longlist}[(A5)]
\item[(A1)] There exists a positive constant $C_1$ such that
    $\lambda_{\min} (n^{-1}\mathbf{X}_{A_0}^T\mathbf{X}_{A_0} )\geq
    C_1$.

\item[(A2)] The random errors $\varepsilon_1,\ldots,\varepsilon_n$ are
    i.i.d. mean zero sub-Gaussian random variables with a scale factor
    $0<\sigma<\infty$, that is, $ E[\exp(t\varepsilon_i)]\leq
    e^{\sigma^2t^2/2}, \forall t. $

\item[(A3)] The penalty function $p_{\lambda}(t)$ is assumed to be
    increasing and concave for $t\in[0,+\infty)$ with a continuous
    derivative $\dot{p}_{\lambda}(t)$ on $(0,+\infty)$. It admits a
    convex-concave decomposition as in (\ref{decom}) with
    $J_{\lambda}(\cdot)$ satisfies: $\nabla J_\lambda(\llvert t\rrvert
    )=-\lambda\operatorname{sign}(t)$ for $|t|>a\lambda$, where $a>1$
    is a constant; and $|\nabla J_\lambda(\llvert t\rrvert )|\leq|t|$
    for $|t|\leq b\lambda$, where $b\leq a$ is a positive constant.

\item[(A4)] The design matrix $\mathbf{X}$ satisfies: $\gamma=
    \min_{\bolds{\delta}\ne0, \|\bolds{\delta}_{A_0^c}\|_1 \leq3
    \|\bolds{\delta}_{A_0}\|_1} \frac{\|\mathbf{X}\bolds{\delta
    }\|}{\sqrt{n} \| \bolds{\delta}_{A_0}\|}>0$.

\item[(A5)] Assume that $\lambda=o(d_*)$ and $\tau=o(1)$, where $d_*$
    is defined on page~\pageref{page5}, $\lambda$~and~$\tau$ are the two parameters in
    the modified CCCP algorithm given in the first paragraph of
    Section~\ref{sec2.3}.
\end{longlist}

\begin{rem}
Condition (A1) concerns the true model and is a common assumption in
the literature on high-dimensional regression. Condition (A2) implies
that for a vector $\mathbf{a}=(a_1,\ldots,a_n)^T$,
%
%
%
%
\begin{equation}
\label{dog1} P \bigl( \bigl|\mathbf{a}^T\bolds{\varepsilon} \bigr|>t \bigr)
\leq2 \exp \biggl(-\frac{t^2}{2\sigma^2\|\mathbf{a}\|^2} \biggr),\qquad t\geq0.
\end{equation}
Condition (A3) is satisfied by popular nonconvex penalty functions such
as SCAD and MCP. Note that the condition $\nabla J_\lambda(\llvert
t\rrvert )=-\lambda\operatorname{sign}(t)$ for $|t|>a\lambda$ is
equivalent to assuming that $\dot{p}_{\lambda}(\llvert t\rrvert )=0$,
$\forall |t|>a\lambda$, that is, large coefficients are not penalized,
which is exactly the motivation for nonconvex penalties. Condition
(A4), which is given in \citet{B2}, ensures a desirable bound on
the $L_1$ estimation loss of the Lasso estimator. Note that the CCCP
algorithm yields the Lasso estimator after the first iteration, so the
asymptotic properties of the CCCP estimator is related to that of the
Lasso estimator. Condition (A4) holds under the restricted eigenvalue
condition which is known to be a relatively mild condition on the
design matrix for high-dimensional estimation. In particular, it is
known to hold in some examples where the covariates are highly
dependent, and is much weaker than the irrepresentable condition
[\citet{Zhao1}] which is almost necessary for Lasso to be model
selection consistent.
\end{rem}

\subsection{Property of the solution path}\label{sec3.1}

We first state a useful lemma that characterizes a nonasymptotic
property of the oracle estimator in high dimension. The result is an
extension of that in \citet{K1} under the more general
sub-Gaussian random error condition.

%
\begin{lemma}\label{lem1} For any given ${0<b_1<1}$ and $0<b_2<1$,
consider the events
\[
F_{n1}= \Bigl\{\max_{j\in A_0}\bigl|\widehat{\beta}{}_j^{(o)}-\beta_j^*\bigr|\leq b_1
\lambda \Bigr\}\quad\mbox{and}\quad F_{n2}= \Bigl\{\max
_{j\in A_0^c}\bigl|S_j \bigl(\widehat{\bolds{\beta}}{}^{(o)} \bigr)\bigr|\leq b_2\lambda \Bigr\},
\]
where
$S_j(\bolds{\beta})=-n^{-1}\mathbf{x}_{(j)}^T(\mathbf{y}-\mathbf
{X}\bolds{\beta})$. Then under conditions \textup{(A1)} and
\textup{(A2)},
\[
P (F_{n1}\cap F_{n2} ) \geq1- 2q{\exp \bigl[-C_1b_1^2n
\lambda^2/ \bigl(2\sigma^2 \bigr) \bigr]}-2(p-q)\exp
\bigl[-nb_2^2\lambda^2/ \bigl(2
\sigma^2 \bigr) \bigr].
\]
\end{lemma}

The proof of Lemma~\ref{lem1} is given in the online supplementary
material [\citet{W3}].

Theorem \ref{main} below provides a nonasymptotic bound of the
probability the solution path contains the oracle estimator. Under
general conditions, this probability tends to one.

%
\begin{theorem}\label{main}
\textup{(1)} Assume that conditions \textup{(A1)--(A5)} hold. If
$\tau\gamma^{-2}q=o(1)$,
%
then for all $n$ sufficiently large,
\[
P \bigl(\widehat{\bolds{\beta}}(\lambda)=\widehat{\bolds{\beta}}{}^{(o)}
\bigr) \geq{1-8p\exp \bigl(-n\tau^2\lambda^2/ \bigl(8
\sigma^2 \bigr) \bigr)}.
\]

\textup{(2)} Assume that conditions \textup{(A1)--(A5)} hold. If
{$n\tau^2\lambda^2\rightarrow\infty$, $\log p=o(n\tau^2\lambda^2)$ and
$\tau\gamma^{-2}q=o(1)$}, then
\[
P \bigl(\widehat{\bolds{\beta}}(\lambda)=\widehat{\bolds{\beta}}{}^{(o)}
\bigr)\rightarrow1
\]
as $n\rightarrow\infty$.
\end{theorem}


\begin{rem}
\citet{Me} considered thresholding Lasso, which has the oracle
property under an incoherent design condition in the ultra-high
dimension. \citet{Zhou2} further proposed and investigated a
multi-step thresholding procedure which can accurately estimate the
sparsity pattern under the restricted eigenvalue condition of
\citet{B2}. These theoretical results are derived by assuming the
initial Lasso is obtained using a theoretical tuning parameter value,
which depends on the unknown random noise variance $\sigma^2$.
Estimating $\sigma^2$ is a difficult problem in high-dimensional
setting, particularly when the random noise is non-Gaussian. On the
other hand, if the true value of $\sigma^2$ is known a priori, then it
is possible to derive variable selection consistency under somewhat
more relaxed conditions on the design matrix than those in the current
paper. Adaptive Lasso, originally proposed by \citet{Z1} for fixed
dimension, was extended to high dimension by \citet{H2} under a
rather strong mutual incoherence condition. \citet{Zhou1} derived
the consistency of adaptive Lasso in high dimension under similar
conditions on $X$, but still requires complex conditions on $s$ and
$d_*$. Some favorable empirical performance of the multi-step
thresholded Lasso versus the adaptive Lasso was reported in
\citet{Zhou2}. A theoretical comparison of these two procedures in
high dimension was considered by \citet{Geer2} and Chapter~7 of
\citet{B3}. For both adaptive and thresholded Lasso, if a
covariate is deleted in the first step, it will be excluded from the
final selected model. \citet{Zh1} proved that selection
consistency holds for the MCP solution at the universal penalty level
$\sigma\sqrt{2\log p/n}$. The LLA algorithm, which \citet{Z2}
originally proposed for fixed dimensional models, alleviates this
problem and has the potential to be extended to the ultra-high
dimension under conditions similar as those in this paper. Needless to
say, the performances of the above procedures all depend on the choice
of tuning parameter. However, the important issue of tuning parameter
selection has not been addressed.
\end{rem}

\begin{rem}
We proved that the calibrated CCCP algorithm which involves merely two
iterations is guaranteed to yield a solution\vadjust{\goodbreak} path that contains the
oracle estimator with high probability under general conditions. To
provide some intuition on this theory, we first note that the first
step of the algorithm yields the Lasso estimator, albeit with a small
penalty level $\tau\lambda$. If we denote the first step estimator by
$\widehat{\beta}{}_j^{(\mathrm{Lasso})}(\tau\lambda)$, then based on
the optimization theory, the oracle property is achieved when
\begin{eqnarray*}
&\displaystyle \min_{j\in A_0}\bigl|\widehat{\beta}{}_j^{(\mathrm{Lasso})}(
\tau\lambda)\bigr|\geq a\lambda>\lambda,&
\\
&\displaystyle \operatorname{sign} \bigl(\widehat{\beta}{}_j^{(o)}
\bigr)= \operatorname{sign} \bigl(\beta_j^* \bigr),\qquad j\in A_0,&
\\
&\displaystyle \max_{j\notin A_0} \bigl\llvert\nabla J_\lambda \bigl(
\widehat{\beta}{}_j^{(\mathrm{Lasso})}(\tau\lambda) \bigr) \bigr
\rrvert+n^{-1}\bigl\|\mathbf{X}^T_{A_0^c}(\mathbf{Y}-
\mathbf{X})\widehat{\bolds{\beta}}{}^{(o)}\bigr\|_{\infty
}\leq\lambda.&
\end{eqnarray*}
The proof of Theorem \ref{main} relies on the following condition:
%
%
%
%
\begin{equation}
\label{key} \bigl\|\widehat{\bolds{\beta}}{}^{(\mathrm{Lasso})}(\tau\lambda)-\bolds {\beta
}^*\bigr\|_{\infty}\leq\lambda/2,\qquad\min_{\beta_j^*\neq0}\bigl|
\beta_j^*\bigr|>a\lambda+\lambda/2
\end{equation}
for the given $a>1$. The proof proceeds by bounding the first part of
(\ref{key}) using a result\vspace*{1pt} of \citet{B2} via
$\|\widehat{\bolds{\beta}}{}^{(\mathrm{Lasso})}(\tau\lambda)-\bolds
{\beta}\|_{\infty}
\leq\|\widehat{\bolds{\beta}}{}^{(\mathrm{Lasso})}(\tau\lambda)-\bolds
{\beta}\|_2$. In Section~\ref{sec3.3}, we considered an alternative
approach using the recent result of \citet{Zh3}, which leads to
weaker requirement on the minimal signal strength under slightly
stronger assumptions on the design matrix. We also noted that Theorem
\ref{main} holds for any $a>1$, although in the numerical studies we
use the familiar $a=3.7$.

How fast the probability that our
estimator is equal to the oracle estimator approaches one depends on the
sparsity level, the magnitude of the smallest signal, the size of
the tuning parameter and the condition of the design matrix.
Corollary~\ref{main3} below confirms that
the path-consistency can hold in ultra-high dimension.
\end{rem}


%
\begin{corollary}\label{main3}
Assume that conditions \textup{(A1)--(A4)} hold. Suppose there are two
positive constants $\gamma_0$ and $K$ such that $\gamma\geq\gamma_0>0$
and $q<K$. If $d_*=O(n^{-c_1})$ for some $c_1\geq0$ and
$p=O(\exp(n^{c_2}))$ for some $c_2>0$, then
\[
P \bigl(\widehat{\bolds{\beta}}(\lambda)=\widehat{\bolds{\beta}}{}^{(o)}
\bigr)\rightarrow1,
\]
provided $\lambda=O(n^{-c_3})$ for some $c_3>c_1$, $\tau^2
n^{1-2c_3-c_2} \rightarrow\infty$ and $\tau=o(1)$.
\end{corollary}

The above corollary indicates that if the true model is very sparse\break
\mbox{(i.e., $q<K$)} and the design matrix behaves well (i.e., $\gamma\geq
\gamma_0>0$), then we can take $\tau$ to be a sequence that converges
to 0 slowly, for example, $\tau=1/\log n$. On the other hand, if one is
concerned that the true model may not be very sparse
($q\rightarrow\infty$) and the design matrix may not behave very well
($\gamma\rightarrow0$), then an alternative choice is to take $\tau
=\lambda$ which works also quite well in practice. The following
corollary establishes that under some general conditions, the choice of
$\tau=\lambda$ yields a consistent solution path under ultra
high-dimensionality.


%
\begin{corollary}\label{main4}
Assume that conditions \textup{(A1)--(A4)} hold. If $q=O(n^{c_1})$ for
some $ c_1\geq0$, $d_*=O(n^{-c_2})$ for some $c_2\geq0$,
$\gamma=O(n^{-c_3})$ for some $c_3\geq0$, $p=O(\exp(n^{c_4}))$ for some
$0<c_4<1$, $\lambda=O(n^{-c_5})$ for some $\max(c_2,
c_1+2c_3)<c_5<(1-c_4)/4$ and $\tau=\lambda$, then
\[
P \bigl(\widehat{\bolds{\beta}}(\lambda)=\widehat{\bolds{\beta}}{}^{(o)}
\bigr)\rightarrow1.
\]
\end{corollary}

\subsection{Property of the high-dimensional BIC}\label{sec3.2}

Theorem \ref{BIC} below establishes the effectiveness of the HBIC
defined in (\ref{hbic}) for selecting the oracle estimator along a
solution path of the calibrated CCCP.

%
\begin{theorem}[(Property of HBIC)]\label{BIC}
Assume that the conditions of Theorem~\ref{main}(2) hold, and {there
exists a positive constant $\kappa$ such that}
%
%
%
%
\begin{equation}
\label{iden} \lim_{n\rightarrow\infty} \min_{A\nsupseteq
A_0, |A|\leq K_n} \bigl
\{n^{-1}\bigl\|(\mathbf{I}_n-\mathbf{P}_A)
\mathbf{X}_{A_0}\bolds{\beta}^*_{A_0}\bigr\|^2 \bigr\}
\geq\kappa,
\end{equation}
where $\mathbf{I}_n$ denotes the $n\times n$ identity matrix and
$\mathbf{P}_A$ denotes the projection matrix onto the linear
space spanned by the columns of $\mathbf{X}_A$. If
$C_n\rightarrow\infty$, $qC_n\log(p)=o(n)$ and $K_n^2\log(p)\log(n)=o(n)$,
then
\[
P(M_{\widehat{\lambda}}=A_0)\rightarrow1
\]
as $n,p\rightarrow\infty$.
\end{theorem}

\begin{rem}
Condition (\ref{iden}) is an asymptotic model identifiability
condition, similar to that in \citet{C2}. This condition states
that if we consider any model which contains at most $K_n$ covariates,
it cannot predict the response variable as well as the true model does
if it is not the true model. To give some intuition of this condition,
as in \citet{C2}, one can show that for $A\nsupseteq A_0$,
\begin{eqnarray*}
n^{-1}\bigl\|(\mathbf{I}_n-\mathbf{P}_A)
\mathbf{X}_{A_0}\bolds{\beta}^*_{A_0}\bigr\|^2 &\geq&
\lambda_{\min} \bigl(n^{-1}\mathbf{X}_{A_0\cup A}^T
\mathbf{X}_{A_0\cup A} \bigr)\bigl\| \bolds{\beta}^*_{A_0\cap A^c}\bigr\|^2
\\
&\geq& \lambda_{\min} \bigl(n^{-1}\mathbf{X}_{A_0\cup A}^T
\mathbf{X}_{A_0\cup A} \bigr)\min_{\beta_j\neq0}{
\beta_j^{*2}}.
\end{eqnarray*}
The theorem confirms that the BIC criterion for shrinkage parameter
selection investigated in \citet{W2}, \citet{W1} and
\citet{Zh4} can be modified and extended to ultra-high
dimensionality. Carefully examining the proof, it is worth noting that
the consistency of the HBIC only requires a consistent solution path
but does not rely on the particular method used to construct the path.
Hence, the proposed HBIC has the potential to be generalized~to other
settings with ultra-high dimensionality.  The sequence $C_n$ should
diverge to $\infty$ slowly, for example, $C_n = \log(\log n)$, which is
used in our numerical studies.
\end{rem}

\subsection{Relaxing the conditions on the minimal signal}\label{sec3.3}

Theorem~\ref{main}, which is the main result of the paper, implies that
the oracle property of the calibrated CCCP estimator requires the
following lower bound on the magnitude of the smallest nonzero
regression coefficient:
%
%
%
%
\begin{equation}
\label{eq:1} d_* \succ\lambda\succ cq\sqrt{\log p/n},
\end{equation}
where $a\succ b$ means $\lim_{n\rightarrow\infty}a/b=\infty$, and $c$
is a constant that depends on the design matrix $\mathbf{X}$ and other
unknown factors such as $\sigma^2$. When the true model dimension $q$
is fixed, the lower bound for $d_*$ is arbitrarily close to the optimal
lower bound $c\sqrt{ \log p/n}$ for nonconvex penalized approaches
[e.g., \citet{Zh1}]. However, when $q$ is diverging, this bound is
suboptimal. In general, there is a tradeoff between the conditions on
$d_*$ and the conditions on the design matrix. Comparing to the results
in the literature, Theorem~\ref{main} imposes weak conditions on the
design matrix and the algorithm we investigate is transparent. In this
section, we will prove that the optimal lower bound of $d_*$ can be
achieved by the calibrated CCCP procedure under a set of slightly
stronger conditions on the design matrix.

Note that the calibrated CCCP estimator depends on
$\widehat{\bolds{\beta}}{}^{(1)}$, which is the Lasso estimator
obtained after the first iteration of the CCCP algorithm. In fact, the
lower bound of $d_*$ is proportional to the $l_\infty$ convergence rate
of $\widehat{\bolds{\beta}}{}^{(1)}$ to~$\bolds{\beta}^*$, and
condition~(A4) only implies that
$\max_{j}|\widehat{\beta}{}^{(1)}_j-\beta^*_j|$ is proportional to
$O_p(q\sqrt{\log p/n}/\tau)$.~If
%
%
%
%
\begin{equation}
\label{eq:linfty} \max_{j}\bigl|\widehat{\beta}
{}^{(1)}_j- \beta^*_j\bigr|=O_p(\sqrt{
\log p/n}/\tau),
\end{equation}
we can show that {$d_*\succ c\sqrt{\log p/n}/\tau$} for any
$\tau=o(1)$, and hence we can achieve almost the optimal lower bound
for $d_*$. Now, the question is under what conditions inequality
(\ref{eq:linfty}) holds. Let $v_{ij}$ be the $(i,j)$ entry of
$\mathbf{X}^T\mathbf{X}$. \citet{L1} derived the convergence rate
(\ref{eq:linfty}) under the condition of mutual coherence:
%
%
%
%
\begin{equation}
\label{eq:coh} \max_{i\ne j} |v_{ij}| > b/q
\end{equation}
for some constant $b>0$. However, the mutual coherence condition would
be too strong for practical purposes
when $q$ is diverging, since it requires that the pairwise correlations
between all possible pairs are sufficiently small. In this subsection,
we give an alternative condition for (\ref{eq:linfty}) based on the
$l_1$ operation norm
of $\mathbf{X}^T\mathbf{X}$.

We replace condition (A4) with the slightly stronger condition
(A4$^{\prime}$) below. We also introduce an additional condition (A6)
based on the matrix $l_1$ operational norm. For a given $m \times m$
matrix $\mathbf{A}$, the $l_1$ operational norm $\|\mathbf{A}\|_1$ is
defined by $\|\mathbf{A}\|_1=\max_{i=1,\ldots,m} \sum_{j=1}^m
|a_{ij}|$, where $a_{ij}$ is the $(i,j)$th entry of $\mathbf{A}$. Let
\begin{eqnarray*}
\zeta_{\max}(m)&=&\max_{|B|\leq m, A_0\subset B} \bigl\| n^{-1}
\mathbf{X}_B^T\mathbf{X} _B
\bigr\|_1,
\\
\zeta_{\min}(m)&=&\max_{|B|\leq m, A_0\subset B} \bigl\|
\bigl(n^{-1} \mathbf{X}_B^T\mathbf{X}
_B \bigr)^{-1}\bigr\|_1.
\end{eqnarray*}

Condition (A4$^{\prime}$): There exist positive constants $\alpha $ and
$\kappa_{\min}$ such that
%
%
%
%
\begin{equation}
\label{eq:ximin} \xi_{\min} \bigl((\alpha+1)q \bigr) \geq
\kappa_{\min}
\end{equation}
and
%
%
%
%
\begin{equation}
\label{eq:ximax} \frac{\xi_{\max}(\alpha q)}{\alpha} \leq\frac
{1}{576}\kappa_{\min}\biggl(1-3\sqrt{\frac{\xi_{\max}(\alpha q)}{\alpha\kappa_{\min}}} \biggr)^2,
\end{equation}
where {$\xi_{\max}(m)=\max_{|B|\leq m, A_0\subset B} \lambda_{\max
}(n^{-1} \mathbf{X}_B^T\mathbf{X}_B)$.}

Condition (A6): Let $u=\alpha+1$. There exist finite positive constants
$\eta_{\max}$ and $\eta_{\min}$ such that
\[
\limsup_{n\rightarrow\infty} \zeta_{\max}(uq)\leq\eta_{\max}< \infty
\]
and
\[
\limsup_{n\rightarrow\infty} \zeta_{\min}(uq)\leq\eta_{\min}< \infty.
\]
%

\begin{rem}
Similar conditions to condition (A4$^{\prime}$) were considered by
\citet{Me} and \citet{B2} for the $l_2$ convergence of the
Lasso estimator. However, (\ref{eq:ximax}) of condition
(A4$^{\prime}$), which essentially assumes that $\xi_{\max}(\alpha
q)/\alpha$ is sufficiently small, is weaker, at least asymptotically,
than the corresponding condition in \citet{Me} and
\citet{B2}, which assumes that $\xi_{\max}(q+\min\{n,p\})$ is
bounded. \citet{Zh3} proved that $|\{j\dvtx\widehat{\beta}_j\ne0\}
\cup A_0| \leq(\alpha+1) q$ under condition (A4$^{\prime}$). In
addition, condition (A4$^{\prime}$) implies condition (A4) [see
\citet{B2}]. Condition (A6) is not too restrictive. Assume the
$\mathbf{x}_i$'s are randomly sampled from a distribution with mean
$\mathbf{0}$ and covariance matrix $\Sigma$. If the $l_1$ operational
norm of $\Sigma$ and $\Sigma^{-1}$ are bounded, then we have
$\zeta_{\max}(uq)\leq\max_{|B|\leq uq, A_0\subset B}\|\Sigma_B\|_1
+o_p(1)$ and $\zeta_{\min}(uq)\leq\max_{|B|\leq uq, A_0\subset B}
\|\Sigma_B^{-1}\|_1 +o_p(1)$ provided that $q$ does not diverge too
fast. Here $\Sigma_B$ is the $|B|\times|B|$ submatrix whose entries
consist of $\sigma_{jl}$, the $(j,l)$th entry of $\Sigma$, for $j\in B$
and $l\in B$. See Proposition~A.1 in the online supplementary
material [\citet{W3}] of this paper. An example of $\Sigma$ satisfying
$\max_{|B|\leq uq, A_0\subset B}\| \Sigma_B\|_1<\infty$ and
$\max_{|B|\leq uq, A_0\subset B}\|\Sigma_B^{-1}\|_1<\infty$ is a block
diagonal matrix where each block is well posed and of finite dimension.
Moreover, condition (A6) is almost necessary for the $l_\infty$
convergence of the Lasso estimator. Suppose that $p$ is small and $d_*$
is large so that all coefficients of the Lasso coefficients are
nonzero. Then,
\[
\widehat{\bolds{\beta}}{}^{(1)}=\widehat{\bolds{\beta}}{}^{ls}+
\tau\lambda \bigl(\mathbf{X} ^T\mathbf{X}/n \bigr)^{-1}
\bolds{\delta},
\]
where $\widehat{\bolds{\beta}}{}^{ls}$ is the\vspace*{-2pt} least square estimator,
and $\bolds{\delta}=(\delta_1,\ldots,\delta_p)$ with
$\delta_j=\operatorname{sign}(\widehat{\beta}{}^{ls}_j)$. Hence, for
the sup norm between
$\widehat{\bolds{\beta}}{}^{(1)}-\widehat{\bolds{\beta}}{}^{ls}$ to be
the order of $\tau\lambda$, the $l_1$ operational norm of
$(\mathbf{X}^T\mathbf{X}/n)^{-1}$ should be bounded.
\end{rem}

%
\begin{theorem}\label{tiger}
Assume that conditions \textup{(A1)--(A3)}, \textup{(A4$^{\prime}$)},
\textup{(A5)} and \textup{(A6)} hold.
\begin{longlist}[(1)]
\item[(1)] If $\tau=o(1)$, then for all $n$ sufficiently large,
\[
P \bigl(\widehat{\bolds{\beta}}(\lambda)=\widehat{\bolds{\beta}}{}^{(o)}
\bigr)\geq{1-8p\exp \bigl[-n\tau^2\lambda^2/ \bigl(8
\sigma^2 \bigr) \bigr]}.
\]

\item[(2)] {If $\tau=o(1)$ and $\log p=o(n\tau^2\lambda^2)$}, then
\[
P \bigl(\widehat{\bolds{\beta}}(\lambda)=\widehat{\bolds{\beta}}{}^{(o)}
\bigr)\rightarrow1
\]
as $n\rightarrow\infty$.

\item[(3)] Assume that the conditions of \textup{(2)} and (\ref{iden})
    hold. Let $\widehat{\lambda}$ be the tuning parameter selected by HBIC. If
    $C_n\rightarrow\infty$, $qC_n\log(p)=o(n)$,
    $K_n^2\log(p)\log(n)=o(n)$, then $
    P(M_{\widehat{\lambda}}=A_0)\rightarrow1$, as
    $n,p\rightarrow\infty$.
\end{longlist}
\end{theorem}

\begin{rem}
We only need $\tau=o(1)$ in Theorem \ref{tiger} for the probability
bound of the calibrated CCCP estimator, while Theorem~\ref{main}
requires $\tau\gamma^{-2} q =o(1)$. Under the conditions of Theorem
\ref{tiger}, the oracle property of $\widehat{\bolds{\beta} }(\lambda)$
holds when
%
%
%
\begin{equation}
\label{eq:3} d_*\succ\lambda\succ\frac{1}{\tau} \sqrt{\log p/n}.
\end{equation}
Since $\tau$ can converge to 0 arbitrarily slowly (e.g., $\tau=1/\log
n$), the lower bound of $d_*$ given by (\ref{eq:3}), $\sqrt{\log
p/n}/\tau$, is almost optimal.
\end{rem}

\section{Numerical results}\label{sec4}
\subsection{Monte Carlo studies}\label{sec4.1}

We now investigate the sparsity recovery and estimation properties of
the proposed estimator via numerical simulations. We \mbox{compare} the
following estimators: the oracle estimator which assumes the
availability of the knowledge of the true underlying model; the Lasso
estimator (implemented using the R package glmnet); the adaptive Lasso
estimator [denoted by ALasso, \citet{Z1}, Section~2.8 of
\citet{B3}], the hard-thresholded\vadjust{\goodbreak} Lasso estimator [denoted by
HLasso, Section~2.8, \citet{B3}], the SCAD estimator from
the original CCCP algorithm without calibration (denoted by SCAD); the
MCP estimator with $a=1.5$ and $3$. For Lasso and SCAD, 5-fold
cross-validation is used to select the tuning parameter; for ALasso,
sequential tuning as described in Chapter~2 of \citet{B3} is
applied. For HLasso, following a referee's suggestion, we first used
$\lambda$ as the tuning parameter to obtain the initial Lasso
estimator, then thresholded the Lasso estimator using thresholding
parameter $\eta=c\lambda$ for some $c>0$ and refitted least squares
regression. We denote the solution path of HLasso by
$\widehat{\bolds{\beta}}{}^{\mathrm{HL}}(\lambda,c\lambda)$, and apply
HBIC to select $\lambda$. We consider $c=2$ and set $C_n=\log\log n$ in
the HBIC as it is found they lead to overall good performance for
HLasso. The MCP estimator is computed using the R package PLUS with the
theoretical optimal tuning parameter value
$\lambda=\sigma\sqrt{(2/n)\log p}$, where the standard deviation
$\sigma$ {is taken} to be known. For the proposed calibrated CCCP
estimator (denoted by New), we take $\tau=1/\log n$ and set
$C_n=\log\log n$ in the HBIC. We observe that
the new estimator performs similarly if we take $\tau=\lambda$. 
In the following, we report simulation results from two examples.
Results of additional simulations can be found in the online
supplemental file.

\begin{exa}\label{exa1}
We generate a random sample $\{y_i, \mathbf{x}_i\}$, $i=1,\ldots, 100$
from the following linear regression model:
\[
y_i=\mathbf{x}_i^T\bolds{\beta}^*+
\varepsilon_i,
\]
where $\bolds{\beta}^*=(3,1.5,0,0,2,\mathbf{0}_{p-5}^T)^T$ with
$\mathbf{0}_k$ denoting a $k$-dimensional vector of zeros, the
$p$-dimensional vector $\mathbf{x}_i$ has the $N(\mathbf{0}_p,
\bolds{\Sigma})$ distribution with covariance matrix $\bolds{\Sigma}$,
$\varepsilon_i$ is independent of $\mathbf{x}_i$ and has a normal
distribution with mean zero and standard deviation $\sigma=2$. This
simulation setup was considered in \citet{F1} for a small $p$
case. In this example, we consider $p=3000$ and the following choices
of $\bolds{\Sigma}$:
(1) Case~1a: the $(i,j)$th entry of
$\bolds{\Sigma}$ is equal to $0.5^{|i-j|}$, $1\leq i, j\leq p$;
(2) Case~1b: the $(i,j)$th entry of $\bolds{\Sigma}$ is equal to
$0.8^{|i-j|}$, $1\leq i, j\leq p$;
(3) Case~1c: the $(i,j)$th entry of
$\bolds{\Sigma}$ equal to 1 if $i=j$ and $0.5$ if $1\leq i\neq j\leq p$.
\end{exa}

\begin{exa}\label{exa2}
We consider a more challenging case by modifying Example~\ref{exa1} case~1a. We
divide the $p$ components of $\bolds{\beta}^*$ into continuous blocks
of size 20. We randomly select 10 blocks and assign each block the
value $(3,1.5,0,0,2,\mathbf{0}_{15}^T)/1.5$. Hence, the number of
nonzero coefficients is 30. The entries in other blocks are set to be
zero. We consider $\sigma=1$. Two different cases are investigated:
(1) Case~2a: $n=200$ and $p=3000$;
(2) Case~2b: $n=300$ and $p=4000$.

In the two examples, based on 100 simulation runs we report the average
number of nonzero coefficients correctly estimated to be nonzero (i.e.,
true positive, denoted by TP) and average number of zero coefficients
incorrectly estimated to be nonzero (i.e., false positive, denoted by
FP) and the proportion of times the true model is exactly identified
(denoted by TM). These three quantities describe the ability of various
estimators for sparsity recovery. To measure the estimation accuracy,
we report the mean squared error (MSE), which is defined to be
$100^{-1}\sum_{m=1}^{100}\| \widehat{\bolds{\beta}}{}^{(m)}-\bolds
{\beta}^*\|^2$, where $\widehat{\bolds{\beta}}{}^{(m)}$ is the
estimator from the $m$th simulation run.

%
\begin{table}
\tabcolsep=0pt
 \caption{Example 1. We report TP (the average number of
nonzero coefficients correctly estimated to be nonzero, i.e., true
positive), FP (average number of zero coefficients incorrectly
estimated to be nonzero, i.e., false positive), TM (the proportion of
the true model being exactly identified) and MSE}\label{table1}
\begin{tabular*}{\tablewidth}{@{\extracolsep{\fill}}lccd{2.2}cc@{}}
\hline
\textbf{Case} & \textbf{Method} & \textbf{TP} & \multicolumn{1}{c}{\textbf{FP}} & \textbf{TM} & \textbf{MSE}\\
\hline
1a
&Oracle &3.00 &0.00 &1.00 &0.146\\
&Lasso &3.00 &28.99 &0.00 &1.101\\
&ALasso &3.00 &11.47 &0.01 &1.327\\
&HLasso &3.00 &0.49 &0.79 &0.383\\
&SCAD &3.00 &10.12 &0.08 &1.496\\
&MCP ($a=1.5$) &2.89 &0.28 &0.76 &0.561\\
&MCP ($a=3$) &2.91 &0.42 &0.68 &1.292\\
&New & 2.99 &\multicolumn{1}{c}{\phantom{0}\textbf{0.09}} &\textbf{0.91} &\textbf{0.222}
\\[3pt]
1b
&Oracle & 3.00 &0.00 &1.00 &0.314\\
&Lasso &3.00 &20.64 &0.00 &1.248\\
&ALasso &3.00 &8.84 &0.02 &1.527\\
&HLasso &2.79 & 0.50 & 0.56 & 1.244 \\
&SCAD &2.99 &7.42 &0.17 &1.598\\
&MCP ($a=1.5$) &2.02 &0.51 &0.06 &5.118\\
&MCP ($a=3$) &1.99 &0.60 &0.02 &5.437\\
&New &2.77 &\multicolumn{1}{c}{\phantom{0}\textbf{0.21}} &\textbf{0.66} &\textbf{1.150}
\\[3pt]
1c
&Oracle &3.00 &0.00 &1.00 &0.195\\
&Lasso &2.99 &28.22 &0.00 &2.987\\
&ALasso &2.96 &10.09 &0.02 &2.433\\
&HLasso &2.84 & 0.77 & 0.56 & 1.361 \\
&SCAD &2.96 &18.09 &0.01 &3.428\\
&MCP ($a=1.5$) &2.67 &\multicolumn{1}{c}{\phantom{0}\textbf{0.17}} &\textbf{0.72} &1.636\\
&MCP ($a=3$) &2.77 &0.22 &0.68 &1.677\\
&New & 2.79 & 0.46 & 0.58 & \textbf{1.244} \\
\hline
\end{tabular*}
\end{table}

%
\begin{table}
\tabcolsep=0pt
\tablewidth=250pt
 \caption{Example 2.
Captions are the same as those in
Table~\protect\ref{table1}}\label{table2}
\begin{tabular*}{\tablewidth}{@{\extracolsep{\fill}}lccd{3.2}cc@{}}
\hline
\textbf{Case} & \textbf{Method} & \textbf{TP} & \multicolumn{1}{c}{\textbf{FP}} & \textbf{TM} & \textbf{MSE} \\
\hline
2a
&Oracle &30.00 &0.00 &1.00 &0.223\\
&Lasso &30.00 &143.14 &0.00 &3.365\\
&ALasso &29.98 &7.50 &0.00 &0.393\\
&HLasso &29.97 & 1.09 & 0.74 & 0.312 \\
&SCAD &29.98 &46.15 &0.00 &2.495\\
&MCP ($a=3$) &29.83 &0.50 &\textbf{0.92} &0.807\\
&New &29.99 & \multicolumn{1}{c}{\phantom{00}\textbf{0.20}} & 0.89 & \textbf{0.247}
\\[3pt]
2b
&Oracle &30.00 & 0.00 &1.00 &0.137\\
&Lasso &30.00 &133.65 &0.00 &1.089\\
&ALasso &30.00 &1.32 &0.29 &0.165\\
&HLasso &30.00 & \multicolumn{1}{c}{\phantom{00}\textbf{0.00}} & \textbf{1.00} & 0.137 \\
&SCAD &30.00 &21.83 &0.00 &0.599\\
&MCP ($a=3$) &30.00 &0.08 &0.92 &0.137\\
&New &30.00 &\multicolumn{1}{c}{\phantom{00}\textbf{0.00}} & 0.99 &\textbf{0.135}\\
\hline
\end{tabular*}
\end{table}

The results are summarized in
Tables \ref{table1}~and~\ref{table2}. It is not surprising that Lasso always overfits.
Other procedures improve the performance of Lasso by reducing the false
positive rate. The SCAD estimator from the original CCCP algorithm
without calibration has no guarantee to find a good local minimum and
has low probability of identifying the true model.
The best overall performance is achieved by the calibrated new
estimator: the probability of identifying the true model is high and
the MSE is relatively small. The HLasso (with thresholding parameter
selected by our proposed HBIC) and MCP (using PLUS algorithm and the
theoretically optimal tuning parameter) also have overall fine
performance. We do not report the results of the MCP with $a=1.5$ for
Example~\ref{exa2} since the PLUS algorithm sometimes runs into
convergence problems.
\end{exa}

\subsection{Real data analysis}\label{sec4.2}
To demonstrate the application, we analyze the gene expression data set
of \citet{S2}, which contains expression values of~31,042 probe
sets on 120 twelve-week-old male offspring of rats. We are interested
in identifying genes whose expressions are related to that of gene
TRIM32 (known to be associated with human diseases of the retina)
corresponding to probe 1389163\_at. We first preprocess the data as
described in \citet{H2} to exclude genes that are either not
expressed or lacking sufficient variation. This leaves 18,957 genes.

%
\begin{table}
\tablewidth=250pt
\tabcolsep=0pt
\caption{Gene expression data analysis. The results are
based on 100 random partitions of the original data set}\label{table3}
\begin{tabular*}{\tablewidth}{@{\extracolsep{\fill}}lcd{2.2}c@{}}\hline
$\bolds{p}$ & \textbf{Method} & \multicolumn{1}{c}{\textbf{ave model size}} & \textbf{Prediction error} \\
\hline
1000
&Lasso &31.17 &\textbf{0.586}\\
&ALasso &11.76 &0.646\\
&HLasso &12.04 &0.676 \\
&SCAD &4.81 &0.827\\
&MCP ($a=1.5$) &11.79 &0.668\\
&MCP ($a=3$) &7.02 &0.768\\
&New &8.50 &0.689
\\[3pt]
2000
&Lasso &32.01 &\textbf{0.604}\\
&ALasso &11.01 &0.661\\
&HLasso &10.82 &0.689 \\
&SCAD &4.57 &0.850\\
&MCP ($a=1.5$) &11.33 &0.700\\
&MCP ($a=3$) &6.78 &0.788\\
&New &7.91 &0.736
\\
\hline
\end{tabular*}
\end{table}

For the analysis, we select 3000 genes that display the largest
variance in expression level. We further analyze the top $p$ ($p=1000$
and 2000) genes that have the largest absolute value of marginal
correlation with gene TRIM32. We randomly partition the 120 rats into
the training data set (80 rates) and testing data set (40 rats). We use
the training data set to fit the model and select the tuning parameter;
and use the testing data set to evaluate the prediction performance. We
perform 1000 random partitions and report in Table~\ref{table3} the average model
sizes and the average prediction error on the testing data set for
$p=1000$ and 2000. For the MCP estimators, the tuning parameters are
selected by cross-validation since the standard deviation of the
random
error is not known. We observe that the Lasso procedure yields the
smallest prediction error. However, this is achieved by fitting
substantially more complex models. The calibrated CCCP algorithm as
well as ALasso and HLasso result in much sparser models with still
small prediction errors. The performance of the MCP procedure is
satisfactory but its optimal performance depends on the parameter $a$.
In screening or diagnostic applications, it is often important to
develop an accurate diagnostic test using as few features as possible
in order to control the cost. The same consideration also matters when
selecting target genes in gene therapies.

We also applied the calibrated CCCP procedure directly to the 18,957
genes and evaluated the predicative performance based on 100 random
partitions. The calibrated CCCP estimator has an average model size 8.1
and an average prediction error 0.58. Note that the model size and
predictive performance are similar to what we obtain when we first
select 1000 (or 2000) genes with the largest variance and marginal
correlation. This demonstrates the stability of the calibrated CCCP
estimator in ultra-high dimension.

When a probe is simultaneously identified by different variable
selection procedures, we consider it as evidence for the strength of
the signal. Probe 1368113\_at is identified by both Lasso and the
calibrated CCCP estimator. This probe corresponds to gene tff2, which
was found to up-regulate cell proliferation in developing mice retina
[\citet{P1}]. On the other hand, the probes identified by the
calibrated CCCP but not by Lasso also merit further investigation. For
instance, probe 1371168\_at was identified by the calibrated CCCP
estimator but not by Lasso. This probe corresponds to gene mpp2, which
was found to be related to protein metabolism abnormalities in the
development of retinopathy in diabetic mice [\citet{Gao}].


\subsection{Extension to penalized logistic regression}\label{sec4.3}
Regularized logistic regression is known to automatically result in a
sparse set of features for classification in ultra-high dimension
[\citet{Geer}, \citet{K3}]. We consider the representative
two-class classification problem, where the response variable $y_i$
takes two possible values 0 or 1, indicating the class membership. It
is assumed that
%
%
%
%
\begin{equation}
\label{binary} 
P(y_i=1\mid
\mathbf{x}_i)=\exp \bigl(\mathbf{x}_i^T
\bolds{\beta} \bigr)/ \bigl\{1+\exp \bigl(\mathbf{x}_i^T
\bolds{\beta} \bigr) \bigr\}.
\end{equation}
The penalized logistic regression estimator minimizes
\[
n^{-1}\sum_{i=1}^n \bigl[-
\bigl(\mathbf{x}_i^T\bolds{\beta} \bigr)y_i+
\log \bigl\{ 1+\exp \bigl(\mathbf{x}_i^T\bolds{\beta}
\bigr) \bigr\} \bigr]+\sum_{j=1}^pp_{\lambda}
\bigl(\llvert \beta_j\rrvert \bigr).
\]
When a nonconvex penalty is adopted, it is easy to see that the CCCP
algorithm can be extended to this case without difficulty as the
penalized log-likelihood naturally possesses the convex-concave
decomposition discussed in Section~\ref{sec2.2} of the main paper,
because of the convexity of the negative log-likelihood for the
exponential family. For easy implementation, the CCCP algorithm can be
combined with the iteratively reweighted least squares algorithm for
ordinary logistic regression, thus taking advantage of the CCCP
algorithm for linear regression. Denote the nonconvex penalized
logistic regression estimator by $\widehat{\bolds{\beta}}$, then for a
new\vspace*{1pt} feature vector $\mathbf{x}$, the predicted class membership is
$I(\exp(\mathbf{x}^T\widehat{\bolds{\beta}})/(1+\exp(\mathbf
{x}^T\widehat{\bolds{\beta}}))>0.5)$.

We demonstrate the performance of nonconvex penalized logistic
regression for classification through the following example: we
generate $\mathbf{x}_i$ as in Example~\ref{exa1} of the main paper, and
the response variable $y_i$ is generated according to (\ref{binary})
with $\bolds{\beta}^*=(3,1.5,0,0,2,\mathbf{0}_{p-50}^T)^T$. We consider
sample size $n=300$ and feature dimension $p=2000$. Furthermore, an
independent test set of size 1000 is used to evaluate the
misclassificaiton error. The simulation results are reported in
Table~\ref{table4}. The results demonstrate that the calibrated CCCP
estimator is effective in both accurate classification and identifying
the relevant features.

We expect that the theory we derived for the linear regression case
continues to hold for the logistic regression under similar conditions
due to the convexity of the negative log-likelihood function and the
fact that the Bernoulli random variables automatically satisfies the
sub-Gaussian tail assumption. The latter is essential for obtaining the
exponential bounds in deriving the theory.

%
\begin{table}
\tabcolsep=0pt
\tablewidth=250pt
\caption{Simulations for classification in high dimension ($n=300$, $p=2000$)}\label{table4}
\begin{tabular*}{\tablewidth}{@{\extracolsep{\fill}}lcd{2.2}cc@{}}
\hline
\textbf{Method} & \textbf{TP} & \multicolumn{1}{c}{\textbf{FP}} & \textbf{TM} & \textbf{Misclassification rate} \\
\hline
Oracle & 3.00 &0.00 &1.00 &0.116\\
Lasso & \textbf{3.00} &46.48 &0.00 &0.134\\
SCAD & 2.08 &4.02 &0.04 &0.161\\
ALASSO & 2.02 &4.58 &0.00 &0.188\\
HLASSO &2.87 &\multicolumn{1}{c}{\phantom{0}\textbf{0.00}} &0.87 &0.120\\
MCP ($a=3$) &2.96 &0.56 &0.54 &0.128\\
New & 2.99 & \multicolumn{1}{c}{\phantom{0}\textbf{0.00}} & \textbf{0.99} & \textbf{0.116} \\
\hline
\end{tabular*}
\end{table}

\section{Revisiting local minima of nonconvex penalized regression}\label{sec5}

In the following, we shall revisit the issue of multiple local minima
of nonconvex penalized regression. We derive an $L_2$ bound of the
distance between a sparse local minimum and the oracle estimator. The
result indicates that a local minimum which is sufficiently sparse
often enjoys fairly accurate estimation even when it is not the oracle
estimator. This result, to our knowledge, is new in the literature on
high-dimensional nonconvex penalized regression.

Our theory applies the necessary condition for the local minimizer as
in \citet{T2} for convex differencing problems. Let
\[
Q_n(\bolds{\beta})=(2n)^{-1}\|\mathbf{y}-\mathbf{X}\bolds{
\beta}\| ^2+\lambda\sum_{j=1}^p
p_\lambda\bigl(\llvert \beta_j\rrvert \bigr)
\]
and
\[
\nabla(\bolds{\beta})= \bigl\{\bolds{\xi}\in\mathcal{R}^{p}\dvtx
\xi_j=-n^{-1}{\mathbf{x}_{(j)}^T}(
\mathbf{y}-\mathbf{X}\bolds{\beta})+\lambda l_j \bigr\},
\]
where $l_j=\operatorname{sign}(\beta_j)$ if $\beta_j\neq0$ and
$l_j\in[-1,1]$ otherwise, $1\leq j\leq p$. As $Q_n(\bolds{\beta})$ can
be expressed as the difference of two convex functions, a necessary
condition for $\bolds{\beta}$ to be a local minimizer of $Q_n(\bolds
{\beta})$ is
%
%
%
%
\begin{equation}
\label{locnec} \frac{\partial h_n(\bolds{\beta
})}{\partial\bolds{\beta}} \in\nabla(\bolds{\beta}),
\end{equation}
where $h_n(\bolds{\beta})=\sum_{j=1}^p J_\lambda(\llvert \beta_j\rrvert
)$, where $J_\lambda(\llvert \beta_j\rrvert )$ is defined in
Section~\ref{sec2.2} for SCAD and MCP penalty functions.

To facilitate our study, we introduce below a new concept.

%
\begin{definition}\label{de5.1}
The relaxed sparse Riesz condition (SRC) in an \mbox{$L_0$-}neigh\-borhood of
the true model is satisfied for a positive integer $m$ ($2q\leq m\leq
n$) if
\[
\xi_{\min}(m)\geq c_* \qquad\mbox{for some } 0<c_*<\infty,
\]
where $\xi_{\min}$ is defined in (\ref{ximin}).
\end{definition}

\begin{rem}
The \textit{relaxed SRC condition} is related to, but generally weaker
than the \textit{sparse Reisz condition} [\citet{Zh2},
\citet{Zh1}], the \textit{restricted eigenvalue condition} of
\citet{B2} and the \textit{partial orthogonality condition} of
\citet{H2}.

The theorem below unveils that for a given sparse estimator which is a
local minimum of (\ref{model2}), its $L_2$ distance to the oracle
estimator $\widehat{\bolds{\beta}}{}^{(o)}$ has an upper bound, which
is determined by three key factors: tuning parameter $\lambda$, the
sparsity size of the local solution, and the magnitude of the smallest
sparse eigenvalue as characterized by the relaxed SRC condition. To
this end, we consider any local minimum
$\widehat{\bolds{\beta}}=(\widehat{\beta}_j,\ldots,\widehat{\beta}_j)^T$
corresponding to the tuning parameter $\lambda$. Assume that the
sparsity size of this local solution satisfies:
$\|\widehat{\bolds{\beta}}\|_0\leq qu_n$ for some $u_n>0$.
\end{rem}


%
\begin{theorem}[(Properties of the local minima of
nonconvex penalized regression)]\label{local}
Consider SCAD or MCP penalized least
squares regression. Assume that conditions \textup{(A1)} and
\textup{(A2)} hold,
and that the relaxed SRC condition in an \mbox{$L_0$-}neighborhood of the true
model is satisfied for $m=qu_n^*$ where $u_n^*=u_n+1$. Then
if $\lambda=o(d_*)$, then for all $n$ sufficiently large,
%
%
%
%
\begin{eqnarray}\label{L2bound}
&& P \Bigl\{\bigl\|\widehat{\bolds{\beta}}(\lambda)-\widehat{\bolds{\beta}}{}^{(o)}\bigr\|\leq2\lambda\sqrt{qu_n^*}
\xi_{\min}^{-1} \bigl(qu_n^* \bigr) \Bigr\}\nonumber
\\
&&\qquad \geq 1- 2q\exp \bigl[-C_1n(d_*-a\lambda)^2/ \bigl(2
\sigma^2 \bigr) \bigr]
\\
&&\quad\qquad{} -2(p-q)\exp \bigl[-n\lambda^2/
\bigl(2 \sigma^2 \bigr) \bigr],\nonumber
\end{eqnarray}
where $\xi_{\min}(m)$ is defined in (\ref{ximin}) and the positive
constant $C_1$ is defined in~\textup{(A1)}.
\end{theorem}

%
\begin{corollary}\label{col_local}
Under the conditions of Theorem \ref{local}, if we take $\lambda
=\sqrt{3\log(p)/n}$, then we have
%
\begin{eqnarray*}
&&P \biggl\{\bigl\|\widehat{\bolds{\beta}}(\lambda)-\widehat{\bolds{\beta
}}{}^{(o)}\bigr\|^2\leq12\frac{qu_n^*\log(p)}{n\xi_{\min}^{2}(qu_n^*)} \biggr\}
\\
&&\qquad \geq 1- 2q\exp \bigl[-C_1n(d_*-a\lambda)^2/ \bigl(2
\sigma^2 \bigr) \bigr]-2(p-q)\exp \bigl[-n\lambda^2/
\bigl(2 \sigma^2 \bigr) \bigr].
\end{eqnarray*}
\end{corollary}

The simple form in the above corollary suggests that if a local minimum
is sufficiently sparse, in the sense that $u_n$ diverge to $\infty$
very slowly, this bound is nevertheless quite tight as the rate
$q\log(p)/n$ is near-oracle. The factor $u_n\xi_{\min}^{-2}(qu_n^*)$ is
expected to go to infinity at a relatively slow rate if the local
solution is sufficiently sparse. Our experience with existing
algorithms for solving nonconvex penalized regression is that they
often yield a sparse local minimum, which however has a low probability
to be the oracle estimator itself.

\section{Proofs}\label{sec6}
We will provide here proofs for the main theoretical results in this
paper.

\begin{pf*}{Proof of Theorem \ref{main}}
By definition,
$\widehat{\bolds{\beta}}(\lambda)=\arg\min_{\bolds{\beta}}
Q_{\lambda}(\bolds{\beta}\mid\widehat{\bolds{\beta}}{}^{(1)}), $ where
$Q_{\lambda}(\bolds{\beta}\mid\widehat{\bolds{\beta
}}{}^{(1)})=(2n)^{-1}\|\mathbf{y}-\mathbf{X}\bolds{\beta}\|^2+
\sum_{j=1}^p \nabla J_\lambda(\llvert
\widehat{\beta}{}_{j}^{(1)}\rrvert ) \beta_j +
\lambda\sum_{j=1}^p|\beta_j|$. Since
$Q_{\lambda}(\bolds{\beta}\mid\widehat{\bolds{\beta}}{}^{(1)})$ is a
convex function of $\bolds{\beta}$, the KKT condition is necessary and
sufficient for characterizing the minimum. To verify that
$\widehat{\bolds{\beta}}{}^{(o)}$ is the minimizer of
$Q_{\lambda}(\bolds{\beta}\mid\widehat{\bolds{\beta}}{}^{(1)})$, it is
sufficient to show that
%
%
%
%
\begin{equation}
\label{eq:one} n^{-1} {{\mathbf{x}}_{(j)}^{T}}
\bigl(\mathbf{y}-\mathbf{X}\widehat{\bolds{\beta}}{}^{(o)} \bigr)+\nabla
J_\lambda \bigl(\bigl| \widehat{\beta} {}^{(1)}_{j}\bigr|
\bigr) + \lambda\operatorname{sign} \bigl(\widehat{\beta} {}^{(o)}_j
\bigr)=0,\qquad j\in A_0
\end{equation}
and
%
%
%
%
\begin{equation}
\label{eq:two} \bigl| n^{-1} {{\mathbf{x}}_{(j)}^{T}}
\bigl(\mathbf{y}-\mathbf{X}\widehat{\bolds{\beta}}{}^{(o)} \bigr)+\nabla
J_\lambda \bigl(\bigl| \widehat{\beta} {}^{(1)}_{j}\bigr|
\bigr) \bigr|\leq\lambda,\qquad j\notin A_0.
\end{equation}

We first verify (\ref{eq:one}). Note that with the initial value
$\mathbf{0}$, we have
$\widehat{\bolds{\beta}}{}^{(1)}=\arg\min_{\bolds{\beta}}\{
(2n)^{-1}\|\mathbf{y}-\mathbf{X} \bolds{\beta}\|^2+
\tau\lambda\|\bolds{\beta}\|_1\}$. Let
$F_{n3}=\{\|\widehat{\bolds{\beta}}{}^{(1)}-\bolds{\beta}^*\|_1\leq
16\tau\lambda\gamma^{-2}q\}$, where $\|\cdot\|_1$ denotes the $L_1$
norm. By modifying the proof of Theorem~7.2 of \citet{B2},
we can show that
under the conditions of the theorem,
%
%
%
%
\begin{equation}
\label{l1bound} P(F_{n3})\geq1-2p\exp \bigl(-n\tau^2
\lambda^2/ \bigl(8\sigma^2 \bigr) \bigr).
\end{equation}
By the assumption of the theorem, on the event $F_{n3}$,
$\|\widehat{\bolds{\beta}}{}^{(1)}-\bolds{\beta}^*\|_1\leq\lambda/2$
for all $n$ sufficiently large. {Furthermore, we consider the event
$F_{n1}$ defined in Lemma \ref{lem1} with $b_1=1/2$. By Lemma
\ref{lem1}, we have $P(\llVert
\widehat{\bolds{\beta}}{}^{(o)}-\bolds{\beta}^*\rrVert _{\infty}\leq
\lambda/2)\geq1-2q\exp[-C_1n\lambda^2/(8\sigma^2)]$. By the assumption
$\lambda=o(d_*)$, for all $n$ sufficiently large, on the event
$F_{n1}\cap F_{n3}$, we have\vspace*{-2pt}
$\operatorname{sign}(\widehat{\beta}{}^{(1)}_j)=\operatorname
{sign}(\widehat{\bolds{\beta}}{}_j^{(o)})$, for $j\in A_0$ and
$\min_{j\in A_0}|\widehat{\beta}{}^{(1)}_j|>a\lambda$.}
Hence,\vspace*{-2pt} by condition~(A3), on the event $F_{n1}\cap F_{n3}$, $\nabla
J_\lambda(\llvert \widehat{\beta}{}^{(1)}_j\rrvert
)=-\lambda\operatorname{sign} (\widehat{\beta}{}^{(1)}_j)=
-\lambda\operatorname{sign}(\widehat{\beta}{}_j^{(o)})$. Furthermore,
$n^{-1}\mathbf{x}^T_{(j)}(\mathbf{y}-\mathbf{X}\widehat{\bolds{\beta}}{}^{(o)})=0$,
for $j\in A_0$, following the definition of the oracle estimator.
Therefore, (\ref{eq:one})~holds with probability at least $1-
{2q\exp[-C_1n\lambda^2/(8\sigma^2)]}- 2p\exp
(-n\tau^2\lambda^2/ (8\sigma^2) )$.

Next, we verify (\ref{eq:two}). On the event $F_{n3}$, we have
$\max_{j\notin A_0}|\widehat{\beta}{}^{(1)}_j|\leq\lambda/2$, for all
$n$ sufficiently large. We consider the event $F_{n2}$ defined in
Lemma~\ref{lem1} with $b_2=1/2$. Lemma \ref{lem1} implies that
$P(F_{n2})\geq1-2(p-q)\exp[-n\lambda^2/(8\sigma^2)]$. On the event
$F_{n2}$ we have $\max_{j\in
A_0^c}|n^{-1}{{\mathbf{x}}_{(j)}^{T}}(\mathbf{y}-\mathbf{X}\widehat{\bolds{\beta}}{}^{(o)})|\leq\lambda/2$.
By condition (A3), on the event $F_{n2}\cap F_{n3}$, (\ref{eq:two})
holds, and this occurs with probability at least
$1-2(p-q)\exp[-n\lambda^2/(8\sigma^2)]-2p\exp
(-n\tau^2\lambda^2/(8\sigma^2) )$.

The above two steps proves (1). The result in (2) follows immediately
from (1).
\end{pf*}

\begin{pf*}{Proof of Corollaries \ref{main3}~and~\ref{main4}}
The proof follows immediately from Theorem \ref{main}.
\end{pf*}

\begin{pf*}{Proof of Theorem \ref{BIC}}
Recall that $M_{\lambda }=\{j\dvtx
\widehat{\bolds{\beta}}_j(\lambda)\neq0\}$. We define the following
three index sets: $\Lambda_{n-}=\{\lambda>0\dvtx \lambda\in\Lambda_{n},
A_0\not\subset M_{\lambda}\}$, $\Lambda_{n0}=\{\lambda>0\dvtx
\lambda\in\Lambda_{n}, A_0=M_{\lambda}\}$, and
$\Lambda_{n+}=\{\lambda>0\dvtx \lambda\in\Lambda_{n}, A_0\subset
M_{\lambda}\mbox{ and } A_0\neq M_{\lambda}\}$. In other words,
$\Lambda_{n-}$, $\Lambda_{n0}$ and $\Lambda_{n+}$ denote the sets of
$\lambda$ values which lead to underfitted, exactly fitted and
overfitted models, respectively. For a given model (or equivalently an
index set) $M$, let
$\mathrm{SSE}_{M}=\inf_{\bolds{\beta}_M\in\mathbb{R}^{|M|}}\|\mathbf
{y}-\mathbf{X}_M\bolds{\beta} _{M}\|^2$. That is, $\mathrm{SSE}_{M}$ is
the sum of squared residuals when the least squares method is used to
estimate model $M$. Also, let
$\widehat{\sigma}{}^2_{M}=n^{-1}\mathrm{SSE}_{M}$. From the definition,
we always have
$\widehat{\sigma}{}^2_{\lambda}\geq\widehat{\sigma}{}^2_{M_{\lambda}}$.\vspace*{1pt}

Consider $\lambda_n$ satisfying the conditions of
Theorem~\ref{main}(2). We have $P(M_{\lambda_n}=A_0)\rightarrow1$. We
will prove that $P (\inf_{{\lambda\in\Lambda_{n-}}}
[\mathrm{HBIC}(\lambda)-\mathrm{HBIC}(\lambda_n) ]>0 )\rightarrow1$ and
$P (\inf_{{\lambda\in\Lambda_{n+}}}
[\mathrm{HBIC}(\lambda)-\mathrm{HBIC}(\lambda_n) ]>0 )\rightarrow1$.

\begin{longlist}[Case II]
\item[Case I.] Consider an arbitrary
    $\lambda\in\Lambda_{n-}$, that is, the model corresponding to
    $M_{\lambda}$ is underfitted:
\begin{eqnarray*}
&& P \Bigl(\inf_{{\lambda\in\Lambda_{n-}}} \bigl[\mathrm{HBIC}(\lambda)-
\mathrm{HBIC}(\lambda_n) \bigr]>0 \Bigr)
\\
&&\qquad = P \Bigl(\inf_{{\lambda\in\Lambda_{n-}}} \bigl[\mathrm{HBIC}(\lambda)-
\mathrm{HBIC}(\lambda_n) \bigr]>0, M_{\lambda_n}=A_0
\Bigr)
\\
&&\quad\qquad{} +P \Bigl(\inf_{{\lambda\in\Lambda_{n-}}} \bigl[\mathrm{HBIC}(\lambda)-
\mathrm{HBIC}(\lambda_n) \bigr]>0, M_{\lambda_n}\neq A_0 \Bigr)
\\
&&\qquad \geq P \biggl(\inf_{{\lambda\in\Lambda_{n-}}} \biggl[\log \bigl(\widehat{
\sigma} {}^2_{M_{\lambda}}/\widehat{\sigma} {}^2_{A_0}
\bigr)+\bigl(\llvert M_{\lambda}|-q\bigr) \frac{C_n\log(p)}{n} \biggr]>0
\biggr)+o(1),
\end{eqnarray*}
where the inequality uses Theorem~\ref{main}(2). Furthermore, we
observe that
\[
\log \biggl(\frac{\widehat{\sigma}{}^2_{M_{\lambda}}}{\widehat{\sigma
}{}^2_{A_0}} \biggr) =\log \biggl(1+\frac{n [\widehat{\sigma
}{}^2_{M_{\lambda
}}-\widehat{\sigma}{}^2_{A_0} ]} {
\bolds{\varepsilon}^T(\mathbf{I}_n-\mathbf{P}_{A_0})\bolds
{\varepsilon}}
\biggr).
\]
Applying the inequality $\log(1+x)\geq\min\{0.5x,\log(2)\}$, $\forall
x>0$, we have
\begin{eqnarray*}
&& P \Bigl(\inf_{{\lambda\in\Lambda_{n-}}} \bigl[\mathrm{HBIC}(\lambda)-
\mathrm{HBIC}(\lambda_n) \bigr]>0 \Bigr)
\\
&&\qquad \geq P \biggl(\min \biggl\{\inf_{{\lambda\in\Lambda_{n-}}}\frac
{n (\widehat{\sigma}{}^2_{M_{\lambda}}-\widehat{\sigma
}{}^2_{A_0} )} {
2\bolds{\varepsilon}^T(\mathbf{I}_n-\mathbf{P}_{A_0})\bolds
{\varepsilon}},
\log(2) \biggr\}-\frac{qC_n\log(p)}{n}>0 \biggr)+o(1).
\end{eqnarray*}

To evaluate $\bolds{\varepsilon}^T(\mathbf{I}_n-\mathbf
{P}_{A_0})\bolds{\varepsilon}$, we apply Corollary~1.3 of
\citet{Mi} with their $A_n=\mathbf{I}_n-\mathbf{P}_{A_0}$,
$B_n=2\sigma^4(n-q)$, $\mu_n=\sigma^2$ and $y_n=(n-q)/(\log n)$, we
have $P (\bolds{\varepsilon}^T(\mathbf{I}_n-\mathbf{P}_{A_0})\bolds
{\varepsilon}\leq2\sigma^2(n-q) )\rightarrow1$ as $n\rightarrow\infty$.
Thus
\begin{eqnarray*}
&& P \Bigl(\inf_{{\lambda\in\Lambda_{n-}}} \bigl[\mathrm{HBIC}(\lambda)-
\mathrm{HBIC}(\lambda_n) \bigr]>0 \Bigr)
\\
&&\qquad \geq P \biggl(\min \biggl\{\frac{\inf_{{\lambda\in\Lambda
_{n-}}}n (\widehat{\sigma}{}^2_{M_{\lambda}}-\widehat{\sigma
}{}^2_{A_0} )} {4(n-q)\sigma^2},\log(2) \biggr\}-
\frac{qC_n\log(p)}{n}>0 \biggr)+o(1).
\end{eqnarray*}
In what follows, we will prove that
$qC_n\log(p)=o(\inf_{{\lambda\in\Lambda_{n-}}}n
(\widehat{\sigma}{}^2_{M_{\lambda}}-\widehat{\sigma}{}^2_{A_0} ))$,
which combining with the assumption $qC_n\log(p)=o(n)$ leads to the
conclusion $P (\inf_{{\lambda\in\Lambda_{n-}}} [\mathrm{HBIC}(\lambda
)-\mathrm{HBIC}(\lambda_n) ]>0 )\rightarrow1$.

We have
\begin{eqnarray*}
&& n \bigl(\widehat{\sigma} {}^2_{M_{\lambda}}-\widehat{\sigma }
{}^2_{M_{T}} \bigr)
\\
&&\qquad = \bolds{\mu}^T(\mathbf{I}_n-\mathbf{P}_{M_{\lambda}})
\bolds{\mu}+2\bolds{\mu}^T(\mathbf{I}_n-\mathbf{P}
_{M_{\lambda}})\bolds{\varepsilon} -\bolds{\varepsilon}^T
\mathbf{P}_{M_{\lambda}}\bolds{\varepsilon}+\bolds{\varepsilon}^T
\mathbf{P} _{A_0}\bolds{\varepsilon}
\\
&&\qquad = I_1+I_2-I_3+I_4,
\end{eqnarray*}
where $\bolds{\mu}=\mathbf{X}\bolds{\beta}^*$, $\mathbf
{P}_{M_{\lambda}}$ is the projection matrix into the space spanned by
the columns of $\mathbf {X}_{M_{\lambda}}$, and the definition of
$I_i$, $i=1,2,3,4$, should be clear from the context. Let
$M_{-}=\{j\dvtx j\notin M_{\lambda}, j\in M_T\}$. Note that $M_{-}$ is
nonempty since $M_{\lambda}$ underfits.

By assumption (\ref{iden}), $|I_1|\geq\kappa n$, for all $n$
sufficiently large. To evaluate $I_2$, we have
\[
I_2=2\sqrt{\bolds{\mu}^T(\mathbf{I}_n-
\mathbf{P}_{M_{\lambda
}})\bolds{\mu}}Z(M_{\lambda
})=2\sqrt{I_1}Z(M_{\lambda}),
\]
where $Z(M_{\lambda})=\mathbf{a}^T_n\bolds{\varepsilon}$ with
$\mathbf{a}_n^T=(\bolds{\mu}^T(\mathbf{I}_n-\mathbf{P}_{M_{\lambda
}})\bolds{\mu})^{-1/2}\bolds{\mu}^T(\mathbf{I}
_n-\mathbf{P}_{M_{\lambda}})$. Note that $\|\mathbf{a}_n\|^2=1$ and
$|\Lambda_{-}|\leq\sum_{t=0}^{K_n}{p \choose
t}\leq\sum_{t=0}^{K_n}p^t=\frac{p^{K_n+1}-1}{p-1}\leq2p^{K_n}$.
Applying the sub-Gaussian tail property in (\ref{dog1}), we have
%
\begin{eqnarray*}
&&P \Bigl(\sup_{\eta\in\Lambda_{n-}}\bigl|Z(M_{\lambda})\bigr|>\sqrt{n/\log(n)}
\Bigr)
\\
&&\qquad \leq 4p^{K_n}\exp \bigl(-n/ \bigl(2\sigma^2\log(n) \bigr) \bigr)
\\
&&\qquad = 4\exp \bigl(K_n\log(p)-n/ \bigl(2\sigma^2\log(n)
\bigr) \bigr)\rightarrow0
\end{eqnarray*}
as $K_n\log(p)\log(n)=o(n)$. Hence, $\sup_{\eta\in\Lambda
_{n-}}|I_2|=o(I_1)$. To evaluate $I_3$, let
$r(\lambda)=\operatorname{Trace}(\mathbf{P}_{M_{\lambda}})$. It follows from
Proposition~3 of \citet{Zh1} that for the sub-Gaussian random
variables $\varepsilon_i$, $\forall t>0$,
%
%
%
%
\begin{eqnarray}\label{dog2}
&& P \biggl\{\frac{\bolds{\varepsilon}^T\mathbf{P}_{M_{\lambda
}}\bolds{\varepsilon}}{r(\lambda
)\sigma^2}\geq\frac{1+t}{[1-2/(e^{t/2}\sqrt{1+t}-1)]^2_{+}} \biggr\}
\nonumber\\[-8pt]\\[-8pt]
&&\qquad \leq \exp \biggl(-\frac{r(\lambda)t}{2} \biggr) (1+t)^{(r(\lambda))/2}.\nonumber
\end{eqnarray}
%
We take $t=n/(2\sigma^2K_n\log(n))-1$ in the above inequality. Then
$t\rightarrow\infty$ by the assumptions of the theorem. Thus for all $n$
sufficiently large,
\begin{eqnarray*}
&&P \biggl(\sup_{\lambda\in\Lambda_{n-}}\bigl|\bolds{\varepsilon}^T
\mathbf{P}_{M_{\lambda }}\bolds{\varepsilon}\bigr| >\frac{n}{\log(n)} \biggr)
\\
&&\qquad \leq  P \biggl(\sup_{\lambda\in\Lambda_{n-}} \biggl|\frac{\bolds
{\varepsilon}
^T\mathbf{P}_{M_{\lambda}}\bolds{\varepsilon}} {
r(\lambda)\sigma^2} \biggr| >
\frac{n}{\sigma^2K_n\log(n)} \biggr)
\\
&&\qquad \leq  P \biggl(\sup_{\lambda\in\Lambda_{n-}} \biggl|\frac{\bolds
{\varepsilon}
^T\mathbf{P}_{M_{\lambda}}\bolds{\varepsilon}} {
r(\lambda)\sigma^2} \biggr| >
\frac{1+t}{[1-2/(e^{t/2}\sqrt{1+t}-1)]^2_{+}} \biggr)
\\
&&\qquad \leq  2p^{K_n}\exp \bigl(-n/ \bigl(8\sigma^2K_n
\log(n) \bigr) \bigr) \bigl(n/ \bigl(2\sigma^2K_n\log(n)
\bigr) \bigr)^{K_n/2}
\\
&&\qquad \leq  2\exp \bigl(K_n\log(p)-n/ \bigl(8\sigma^2K_n
\log(n) \bigr)+K_n\log \bigl(n/ \bigl(2\sigma^2K_n
\log(n) \bigr) \bigr) \bigr)
\\
&&\qquad \rightarrow 0,
\end{eqnarray*}
since $K_n^2\log(p)\log(n)=o(n)$. Finally,
$\bolds{\varepsilon}^T\mathbf{P}_{A_0}\bolds{\varepsilon}$ does not
depend on $\lambda$. Similarly as above, $P
(\sup_{\lambda\in\Lambda_{n-}}|I_4|\geq n/\log(n) )\rightarrow0$ by the
sub-Gaussian tail condition. Therefore, with probability approaching
one, $n (\widehat{\sigma}{}^2_{M_{\lambda}}-\widehat{\sigma}{}^2_{A_0}
)$ is dominated by $I_1$. This finishes the proof for the first case as
$qC_n\log(p)=o(n)$.

\item[Case II.] Consider an arbitrary
    $\lambda\in\Lambda_{n+}$, that is, the model corresponding to
    $M_{\lambda}$ is overfitted. In this case, we have
    $\mathbf{y}^T(\mathbf{I}_n-\mathbf{P}_{M_{\lambda}})\mathbf
    {y}=\bolds{\varepsilon}^T(\mathbf{I}_n-\mathbf{P}
    _{M_{\lambda}})\bolds{\varepsilon}$. Therefore, $ n
    (\widehat{\sigma}{}^2_{A_0}-\widehat{\sigma}{}^2_{M_{\lambda}} ) =
    \bolds{\varepsilon}^T(\mathbf{P}_{M_{\lambda}}-\mathbf
    {P}_{A_0})\bolds{\varepsilon}. $
Let $\widehat{\bolds{\varepsilon}}=(\mathbf{I}_n-\mathbf
{P}_{A_0})\bolds{\varepsilon}$, then
\[
\log \biggl(\frac{\widehat{\sigma}{}^2_{A_0}}{\widehat{\sigma
}{}^2_{M_{\lambda}}} \biggr) = \log \biggl(1+\frac{\bolds{\varepsilon
}^T(\mathbf{P}_{M_{\lambda
}}-\mathbf{P}
_{A_0})\bolds{\varepsilon}} {
\bolds{\varepsilon}^T(\mathbf{I}_n-\mathbf{P}_{M_{\lambda
}})\bolds{\varepsilon}}
\biggr) 
\leq\frac{\bolds{\varepsilon}^T(\mathbf{P}_{M_{\lambda}}-\mathbf
{P}_{A_0})\bolds{\varepsilon}} {
\widehat{\bolds{\varepsilon}}{}^T\widehat{\bolds{\varepsilon
}}-\bolds{\varepsilon}^T(\mathbf{P}
_{M_{\lambda}}-\mathbf{P}_{A_0})\bolds{\varepsilon}}
\]
by the fact $\log(1+x)\leq x$, $\forall x\geq0$.

Similarly as in case~I,
\begin{eqnarray*}
&&P \Bigl(\inf_{{\lambda\in\Lambda_{n+}}} \bigl[\mathrm{HBIC}(\lambda)-
\mathrm{HBIC}(\lambda_n) \bigr]>0 \Bigr)
\\
&&\qquad =P \biggl( \inf_{{\lambda\in\Lambda_{n+}}} \biggl[-\log \biggl(
\frac{\widehat{\sigma}{}^2_{A_0}}{\widehat{\sigma
}{}^2_{M_{\lambda}}} \biggr)+\bigl(\llvert M_{\lambda}|-q\bigr)
\frac{C_n\log(p)}{n} \biggr]>0 \biggr)+o(1)
\\
&&\qquad \geq P \biggl(\inf_{{\lambda\in\Lambda_{n+}}} \biggl[\bigl(\llvert
M_{\lambda}|-q\bigr) \frac{C_n\log(p)}{n}-\frac{\bolds
{\varepsilon}^T(\mathbf{P}
_{M_{\lambda}}-\mathbf{P}_{A_0})\bolds{\varepsilon}}{\widehat{\bolds{\varepsilon}}{}^T\widehat{\bolds{\varepsilon}}
-\bolds{\varepsilon}^T(\mathbf{P}_{M_{\lambda}}-\mathbf
{P}_{A_0})\bolds{\varepsilon}} \biggr]>0 \biggr)
\\
&&\quad\qquad{} +o(1)
\\
&&\qquad = P \biggl(\inf_{{\lambda\in\Lambda_{n+}}} \biggl\{\bigl(\llvert
M_{\lambda}|-q\bigr) \biggl[ \frac{C_n\log(p)}{n}-\frac{\bolds{\varepsilon}^T(\mathbf
{P}_{M_{\lambda}}-\mathbf{P}
_{A_0})\bolds{\varepsilon}/(\llvert M_{\lambda}|-q)} {
\widehat{\bolds{\varepsilon}}{}^T\widehat{\bolds{\varepsilon
}}-\bolds{\varepsilon}^T(\mathbf{P}
_{M_{\lambda}}-\mathbf{P}_{A_0})\bolds{\varepsilon}} \biggr]
\biggr\} \biggr)
\\
&&\quad\qquad{} +o(1).
\end{eqnarray*}
It suffices to show that
\[
P \biggl(\inf_{{\lambda\in\Lambda_{n+}}} \biggl[\frac{C_n\log
(p)}{n}-
\frac{\bolds{\varepsilon}^T(\mathbf{P}_{M_{\lambda
}}-\mathbf{P}_{A_0})\bolds{\varepsilon}
/(\llvert M_{\lambda}|-q)}{\widehat{\bolds{\varepsilon}}{}^T\widehat{\bolds{\varepsilon}
}-\bolds{\varepsilon}^T(\mathbf{P}_{M_{\lambda}}-\mathbf
{P}_{A_0})\bolds{\varepsilon}} \biggr]>0 \biggr)\rightarrow1,
\]
which is implied by
\[
P \biggl(\frac{C_n\log(p)}{n}-\frac{\sup_{{\lambda\in\Lambda
_{n+}}}\bolds{\varepsilon}^T(\mathbf{P}_{M_{\lambda}}-\mathbf
{P}_{A_0})\bolds{\varepsilon}
/(\llvert M_{\lambda}|-q)}{\widehat{\bolds{\varepsilon}}{}^T\widehat{\bolds
{\varepsilon}}-\sup_{{\lambda\in\Lambda_{n+}}}\bolds{\varepsilon
}^T(\mathbf
{P}_{M_{\lambda}}-\mathbf{P}
_{A_0})\bolds{\varepsilon}}>0 \biggr)\rightarrow1.
\]
Note that
$E(\widehat{\bolds{\varepsilon}}{}^T\widehat{\bolds{\varepsilon}})=
\operatorname{Var}(\varepsilon_i)\operatorname{Trace}(\mathbf{I}_n-\mathbf{P}_{A_0})\leq
(n-q)\sigma^2$, hence
$\widehat{\bolds{\varepsilon}}{}^T\widehat{\bolds{\varepsilon
}}=O_p(n)$. Similarly as in case~I, we can show that $ P
(\sup_{\lambda\in\Lambda_{n+}}\bolds{\varepsilon}^T(\mathbf
{P}_{M_{\lambda }}-\mathbf{P}_{A_0})\bolds{\varepsilon}
> n/ \log(n) )\rightarrow0$, since $K_n^2\log(p)\log(n)=o(n)$.
Thus, $\widehat{\bolds{\varepsilon}}{}^T\widehat{\bolds{\varepsilon
}}-\sup_{{\lambda\in
\Lambda_{n+}}}\bolds{\varepsilon}^T(\mathbf{P}_{M_{\lambda
}}-\mathbf{P}_{A_0})\bolds{\varepsilon} =O_p(n)$. Furthermore, applying
(\ref{dog2}) by letting $t=8\log(p)-1$, we have for all $n$
sufficiently large,
\begin{eqnarray*}
&&P \biggl(\sup_{\lambda
\in \Lambda_{n+}} \frac{\bolds{\varepsilon}^T(\mathbf{P}_{M_{\lambda
}}-\mathbf{P}_{A_0})\bolds{\varepsilon}}{|M_{\lambda}|-q}> 16 \sigma^2\log(p) \biggr)
\\
&&\qquad \leq \sum_{|M_{\lambda}|=q+1}^{p}{p-q
\choose|M_{\lambda}|-q} \exp \biggl(-\frac{(\llvert M_{\lambda}|-q)t}{2} \biggr)
(1+t)^{(|M_{\lambda}|-q)/2}
\\
&&\qquad = \sum_{k=1}^{p-q}{p-q\choose k} \exp
\bigl(-2k\log(p) \bigr) \bigl(8\log(p) \bigr)^{k/2}
\\
&&\qquad = \sum_{k=1}^{p-q}{p-q \choose k} \biggl(
\frac{\sqrt{8\log(p)}}{p_n^2} \biggr)^k \leq \biggl(1+\frac{\sqrt{8\log(p)}}{p^2}
\biggr)^{p-q}-1\rightarrow0.
\end{eqnarray*}
Thus with
probability approaching one, for all $n$ sufficiently large,
\begin{eqnarray*}
&&\frac{C_n\log(p)}{n}-\frac{\sup_{{\lambda\in
\Lambda_{n+}}}\bolds{\varepsilon}^T(\mathbf{P}_{M_{\lambda
}}-\mathbf{P}_{A_0})\bolds{\varepsilon}
/(\llvert M_{\lambda}|-q)} {
\widehat{\bolds{\varepsilon}}{}^T\widehat{\bolds{\varepsilon
}}-\sup_{{\lambda\in
\Lambda_{n+}}}\bolds{\varepsilon}^T(\mathbf{P}_{M_{\lambda
}}-\mathbf{P}_{A_0})\bolds{\varepsilon}}
\\
&&\qquad > n^{-1}C_n\log(p)-n^{-1}O \bigl(\log(p)\bigr)>0,
\end{eqnarray*}
since $C_n\rightarrow\infty$. This finishes the proof.
\end{longlist}\upqed
\end{pf*}

\begin{pf*}{Proof of Theorem \ref{tiger}}
We will first prove that there exists a constant \mbox{$C>0$} such that for
$F_{n4}=\{\max_{j} |\widehat{\beta}{}^{(1)}_j-\beta^*_j|\leq C
\tau\lambda\}$, we have
%
%
%
%
\begin{equation}\label{F4}
P(F_{n4})\geq1- 2 p\exp \biggl(\frac{-
n\tau^2\lambda^2}{8\sigma^2} \biggr).
\end{equation}
Let $F_{n5}=\{|S_j(\bolds{\beta}^*)|\leq\tau\lambda/2 \mbox{ for all }
j\}$. Since
\[
P \bigl(F_{n5}^c \bigr)\leq\sum
_{j=1}^p P \bigl(\bigl| {\mathbf{x}}_{(j)}^{T}
\bolds{\varepsilon}/n\bigr|> \tau\lambda/2 \bigr)\leq2p \exp \biggl(\frac{-
n\tau^2\lambda^2}{8\sigma^2}
\biggr),
\]
we have
\[
P(F_{n5})\geq1- 2p \exp \biggl(\frac{- n\tau^2\lambda^2}{8\sigma
^2} \biggr).
\]
Hence to prove (\ref{F4}), it suffices to show that $F_{n5}\subset F_{n4}$.

Let
\[
\theta=\inf \biggl\{ \frac{q \|\mathbf{X}^T\mathbf{X}\mathbf{u}\|
_{\infty}}{n\|\mathbf{u}\|
_1}\dvtx \|\mathbf{u}_{A_0^c}
\|_1\leq3 \|\mathbf{u}_{A_0}\|_1 \biggr\}.
\]
Corollary~2 of \citet{Zh3} proves that on the event $F_{n5}$,
$|{A}\cup A_0|\leq(\alpha+1) q$, where $A=\{j\dvtx\widehat{\beta}{}^{(1)}_j\ne0\}$, provided
\[
\frac{\xi_{\max}(\alpha q)}{\alpha} \leq\frac{1}{36} \theta.
\]
Since $\theta\geq\gamma^2/16$ [see (7) of \citet{Zh3}], where
$\gamma$ is defined in (A4) and
\[
\gamma\geq\sqrt{\kappa_{\min}} \biggl(1-3\sqrt{\frac{\xi_{\max
}(\alpha q)}{\alpha\kappa_{\min}}}
\biggr)
\]
[see \citet{B2}], condition (A4$^{\prime}$) implies that
%
%
%
%
\begin{equation}\label{eq:subset}
F_{n5}\subset \bigl\{|{A}\cup
A_0| \leq( \alpha+1) q \bigr\}.
\end{equation}

Let $C(\bolds{\beta})=(2n)^{-1}\|\mathbf{y}-\mathbf{X}\bolds{\beta
}\|^2+ \tau\lambda\sum_{j=1}^p |\beta_j|$.
Then we have
\begin{eqnarray*}
C(\bolds{\beta})-C \bigl(\bolds{\beta}^* \bigr) &=& \sum
_{j=1}^p \bigl(\beta_j-
\beta^*_j \bigr) S_j \bigl(\bolds{\beta}^* \bigr) +
\bigl( \bolds{\beta}-\bolds{\beta}^* \bigr)^{T} \mathbf{X}^{T}
\mathbf{X} \bigl(\bolds{\beta}-\bolds{\beta}^* \bigr)/(2n)
\\
&&{}+ \tau\lambda\sum_{j=1}^p \bigl(|
\beta_j|-\bigl|\beta^*_j\bigr| \bigr).
\end{eqnarray*}
Let $\widehat{\mathbf{X}\bolds{\beta}^*}$ be the projection of
$\mathbf{X}\bolds{\beta}^*$ onto $\operatorname{span}(\mathbf{X}_A)$, the
linear subspace spanned by the column\vspace*{1pt} vectors of $\mathbf{X}_A$. We
define the $p$-dimensional vector $\bolds{\gamma}*$ such that
$\widehat{\mathbf{X}\bolds{\beta}^*}=\mathbf{X}_{A}\bolds{\gamma}^*_A$
and {$\gamma^*_j=0$} for $j\in A^c$. We have
\begin{eqnarray*}
&& \bigl(\widehat{\bolds{\beta}}{}^{(1)}-\bolds{\beta}^*
\bigr)^T \mathbf{X}^T\mathbf{X} \bigl(\widehat{\bolds{
\beta}}{}^{(1)}- \bolds{\beta}^* \bigr)
\\
&&\qquad  = \bigl(\widehat{\bolds{\beta}}{}^{(1)}_A- \bolds{\gamma}_A^*
\bigr)^T \mathbf{X}_A^T \mathbf{X}_A
\bigl(\widehat{\bolds{\beta}}{}^{(1)}_A- \bolds{
\gamma}_A^* \bigr) + \bigl\| \mathbf{X}\bolds{\beta}^*-
\mathbf{X}_A \bolds{\gamma}_A^{*}\bigr\|^2.
\end{eqnarray*}
Therefore, we can write
\begin{eqnarray*}
\widehat{\bolds{\beta}}{}^{(1)} &=&\mathop{\arg\min}_{\bolds{\beta}\dvtx \bolds
{\beta}_{A^c}=\mathbf{0}} \biggl\{
\sum_{j\in A} \beta_j S_j
\bigl(\bolds{\beta}^* \bigr)
\\
&&\hspace*{41pt}{}+ \bigl(\bolds{\beta}_A-\bolds{
\gamma}_A^* \bigr)^T \mathbf{X}^T_A
\mathbf{X}_A \bigl(\bolds{\beta}_A-\bolds{
\gamma}_A^* \bigr)/2n
+\tau\lambda\sum_{j\in A} |\beta_j|\biggr\}.
\end{eqnarray*}
Hence $\widehat{\bolds{\beta}}{}^{(1)}_A -\bolds{\gamma}_A^* = (
\mathbf{X}^T_A\mathbf{X}_A/n )^{-1} \bolds{\theta}_A$, where
$\bolds{\theta}\in R^p$ such that $\theta_j=0$ for $j\in A^c$ and
$\theta_j=-S_j(\bolds{\beta}){-}{\operatorname{sign}(\widehat{\beta}_j)}
\tau\lambda$ for $j\in A$. On $F_{n5}$, $\max_j
|\theta_j|\leq3\tau\lambda/2$. Therefore, condition (A6) with
(\ref{eq:subset}) implies that on the event $F_{n5}$,
%
%
%
%
\begin{equation}
\label{eq:A} \max_{j\in A} \bigl|\widehat{\beta} {}^{(1)}_j-
\gamma_j^*\bigr|\leq\eta_{\min} 3\tau\lambda/2.
\end{equation}
It follows from (\ref{eq:A}) that inequality (\ref{F4}) holds if we
show that $A_0 \subset A$, in which case $\bolds{\gamma}_A^*=\bolds
{\beta}^*_{A}$. We
will prove this by contradiction. Assume $A^{(-)}=A_0\cap A^c$ is
nonempty. Let {$\widehat{\mathbf{x}}_{(j)}$} be the projection of
{$\mathbf{x}_{(j)}$} onto $\operatorname{span}(\mathbf{X}_A)$ and let
{$\tilde{\mathbf{x}}_{(j)}=\mathbf{x}_{(j)}-\widehat{\mathbf{x}}_{(j)}$},
$j\in A^{(-)}$. Then, we can write
\[
\mathbf{X}\bolds{\beta}^*=\mathbf{X}_A\gamma^*_A +
\sum_{j\in
A^-} {\tilde{\mathbf{x}}_{(j)}}
\beta_j^*.
\]
Let $\tilde{\mathbf{y}}=\sum_{j\in A^-} {\tilde{\mathbf{x}}_{(j)}}
\beta_j^*$.
By Lemma \ref{le:low} below, there exists $l\in A^-$ such that
%
%
%
%
\begin{equation}
\label{eq:low} \bigl|{\mathbf{x}_{(l)}^T}\tilde{\mathbf{y}}/n\bigr|
\geq\kappa_{\min} d_*.
\end{equation}
By\vspace*{-2pt} the KKT condition, we have $
\llvert{\mathbf{x}_{(l)}^T}(\mathbf{X}\bolds{\beta}^*-\mathbf
{X}\widehat{\bolds{\beta}}{}^{(1)})/n +S_l(\bolds{\beta}^*)\rrvert
\leq\tau\lambda. $ However we can write
${\mathbf{x}_{(l)}^T}(\mathbf{X}\bolds{\beta}^*-\mathbf{X}\widehat{\bolds{\beta}}{}^{(1)})/n
={\mathbf{x}_{(l)}^T}\mathbf{X}_A(\bolds{\gamma}_A^*-\widehat
{\bolds{\beta}}{}_A^{(1)})/n +
{\mathbf{x}_{(l)}^T}\tilde{\mathbf{y}}/n. $ The inequalities
(\ref{eq:low}) and (\ref{eq:A}) with condition {(A6)} imply that {on
$F_{n5}$}
\begin{eqnarray*}
&& \bigl\llvert{\mathbf{x}_{(l)}^T} \bigl(\mathbf{X}
\bolds{\beta}^*-\mathbf{X}\widehat{\bolds{\beta}}{}^{(1)} \bigr)/n
+S_l \bigl(\bolds{\beta}^* \bigr) \bigr\rrvert
\\
&&\qquad \geq \bigl|{\mathbf{x}_{(l)}^T}\tilde{\mathbf{y}}/n\bigr|-\bigl|{
\mathbf {x}_{(l)}^T}\mathbf{X} _A \bigl(\bolds{
\gamma}_A^*-\widehat{\bolds{\beta}}{}_A^{(1)}
\bigr)/n\bigr| -\bigl|S_l \bigl(\bolds{\beta}^* \bigr)\bigr|
\\
&&\qquad \geq \bigl|\mathbf{x}_{(l)}^T\tilde{\mathbf{y}}/n\bigr|-\bigl\|\mathbf
{X}_{A\cup A_0}^T\mathbf{X}_{A\cup
A_0}\bigr\|_{1}\bigl\|
\bolds{\gamma}_A^*-\widehat{\bolds{\beta}}{}_A^{(1)}
\bigr\|_{\infty} -\bigl|S_l \bigl(\bolds{\beta}^* \bigr)\bigr|
\\
&&\qquad \geq \kappa_{\min} d_* - {\eta_{\max} \eta_{\min}} 3
\tau\lambda/2-\tau\lambda/2 > \tau\lambda
\end{eqnarray*}
if $d_*> 3\tau\lambda({\eta_{\max} \eta_{\min}}+1)/(2\kappa _{\min})$,
which contradicts the KKT condition. Hence, we eventually have
$A_0\subset A$ on $F_{n5}$ and this proves (\ref{F4}).

We now slightly modify the proof of (1) of Theorem~\ref{main}. More
specifically, replacing $F_{n3}$ by $F_{n4}$, we can show that
$F_{n1}\cap F_{n2}\cap F_{n4} \subset\{\widehat{\bolds{\beta}}(\lambda
)=\widehat{\bolds{\beta}}{}^{(o)}\}$, and this proves (1). The result
in (2) follows immediately from (1). The proof of (3) can be done
similarly to that of Theorem \ref{BIC}.
\end{pf*}

In the proof of Theorem~\ref{tiger}, we have used the following lemma,
whose proof is given in the online supplementary material [\citet{W3}].

\begin{lemma}\label{le:low}
There exists $l\in A^-$ which satisfies~(\ref{eq:low}).
\end{lemma}

\begin{pf*}{Proof of Theorem \ref{local}}
By (\ref{locnec}), a local minimizer $\bolds{\beta}$ necessarily
satisfies:
%
%
%
%
\begin{equation}
\label{nce1} -n^{-1}\mathbf{x}_{(j)}^T(\mathbf
{y}-\mathbf{X}\bolds{\beta})+\xi_j=0,\qquad j=1,\ldots,p,
\end{equation}
where $\xi_j=\lambda l_j-\frac{\partial
h_n(\bolds{\beta})}{\partial\beta_j}$, with $l_j=\operatorname
{sign}(\beta_j)$ if $\beta_j\neq0$ and $l_j\in[-1,1]$ otherwise, $1\leq
j\leq p$. It is easy to see that $|\xi_j|\leq\lambda$, $1\leq j \leq
p$. Although the objective function is nonconvex, abusing the notation
a little, we refer to the collection of all vectors in the form of the
left-hand side of (\ref{nce1}) as the subdifferential $\partial
Q_n(\bolds {\beta})$ and refer to a specific element of this set a
subgradient. Then the necessary condition stated above can be
considered as an extension of the classical KKT condition.

Alternatively, minimizing $Q_n(\bolds{\beta})$ can be expressed as a
constrained smooth minimization problem [e.g., \citet{K1}]. By the
corresponding second-order sufficiency of KKT condition [e.g.,
\citet{B1}, page~320], $\widehat{\bolds{\beta}}$ is a~local
minimizer of $Q_n(\bolds{\beta})$ if
\begin{eqnarray*}
n^{-1}{\mathbf{x}_{(j)}^T}(\mathbf{y}-\mathbf{X}
\widehat{\bolds{\beta}})&=&\operatorname{sgn}(\widehat{\beta}_j)
\dot{p}_{\lambda}(\widehat{\beta}_j),\qquad\widehat{\beta}_j\neq0,
\\
n^{-1}\bigl|{\mathbf{x}_{(j)}^T}(\mathbf{y}-\mathbf{X}
\widehat{\bolds{\beta}})\bigr|&\leq& \lambda, \qquad\widehat{\beta}_j= 0.
\end{eqnarray*}
Consider the event $F_n=F_{n2}\cap F_{n6}$, where $F_{n2}$ is defined
in Lemma \ref{lem1} with $b_2=1$, and $F_{n6}= \{\min_{j\in
A_0}|\widehat{\beta}{}_j^{(o)}|\geq a\lambda\}$. Since\vspace*{1.5pt}
$|\widehat{\beta}{}_j^{(o)}| \geq|\beta
_j^*|-|\widehat{\beta}{}_j^{(o)}-\beta_j^*|$ and $\lambda=o(d_*)$,
similarly as in the proof for Lemma \ref{lem1}, we can show that for
all $n$ sufficiently large, $ P(F_{n6})\geq
1-2q\exp[-C_1n(d_*-a\lambda)^2/(2\sigma^2)]. $ By Lemma \ref{lem1}, for
all\vspace*{2pt} $n$ sufficiently large, $P(F_n)\geq1-
2q\exp[-C_1n(d_*-a\lambda)^2/(2\sigma^2)]-2(p-q)\exp[-n\lambda
^2/(2\sigma^2)]$. It is apparent that on the event $F_n$, the oracle
estimator $\widehat{\bolds{\beta}}{}^{(o)}$ satisfies the above
sufficient condition. Therefore, by (\ref{nce1}), there exist
$|\xi_j^{(o)}|\leq\lambda$, $1\leq j\leq p$, such that
\[
-n^{-1}{\mathbf{x}_{(j)}^T} \bigl(\mathbf{y}-
\mathbf{X}\widehat{\bolds{\beta}}{}^{(o)} \bigr)+\xi_j^{(o)}=0.
\]
Abusing notation a little, we denote this zero vector by
$\frac{\partial}{\partial\bolds{\beta}}Q_n(\widehat{\bolds{\beta}}{}^{(o)})$.

Now for any local minimizer $\widehat{\bolds{\beta}}$ which satisfies
the sparsity constraint $\|\widehat{\bolds{\beta}}\|_0\leq qu_n$, we
will prove by contradiction that under the conditions of the theorem we
must have
$\|\widehat{\bolds{\beta}}-\widehat{\bolds{\beta}}{}^{(o)}\|\leq
2\lambda\sqrt{qu_n^*}\xi_{\min}^{-1}(qu_n^*)$, where $u_n^*=u_n+1$.
More specifically, we will derive a contradiction by showing that none of
the subgradients of $Q_n(\bolds{\beta})$ can be zero at
$\bolds{\beta}=\widehat{\bolds{\beta}}$.

Assume instead that
$\|\widehat{\bolds{\beta}}-\widehat{\bolds{\beta}}{}^{(o)}\|>
2\lambda\sqrt{qu_n^*}\xi_{\min}^{-1}(qu_n^*)$. Let $A^*=\{j\dvtx
\widehat{\beta}_j\neq0$ or $\widehat{\beta}{}_j^{(o)}\neq 0\}$,
then $\|\widehat{\bolds{\beta}}_{A^*}-\widehat{\bolds{\beta
}}{}^{(o)}_{A^*}\|>2\lambda\sqrt{qu_n^*}\xi_{\min}^{-1}(qu_n^*)$. Let
$\frac{\partial}{\partial\bolds{\beta}}Q_n(\widehat{\bolds{\beta}})=\break
-n^{-1}{\mathbf{x}_{(j)}^T}(\mathbf{y}-\mathbf{X}\widehat{\bolds
{\beta}})+\eta_j$ be an arbitrary subgradient in the subdifferential
$\partial Q_n(\widehat{\bolds{\beta}})$. Let
$\bolds{\eta}=(\eta_1,\ldots,\eta_p)^T$, then $\eta_j$ satisfies
$|\eta_j|\leq\lambda$, $1\leq j\leq p$. We use
$\frac{\partial}{\partial\bolds{\beta}_{A^*}}Q_n(\widehat{\bolds
{\beta}})$ to denote the size-$|A^*|$ subvector of
$\frac{\partial}{\partial\bolds{\beta}}Q_n(\widehat{\bolds{\beta}})$,
that is, $\frac{\partial}{\partial
\bolds{\beta}_{A^*}}Q_n(\widehat{\bolds{\beta}})= (\frac{\partial
}{\partial\beta_j}Q_n(\widehat{\bolds{\beta}})\dvtx j\in A^*)^T$. And
$\frac{\partial
}{\partial\bolds{\beta}_{A^*}}Q_n(\widehat{\bolds{\beta}}{}^{(o)})$ is
defined similarly. We have
\begin{eqnarray*}
&& \biggl\llvert \biggl(\frac{\partial}{\partial
\bolds{\beta}_{A^*}}Q_n(\widehat{\bolds{\beta}}) \biggr)^T\frac
{(\widehat{\bolds{\beta}}_{A^*}-\widehat{\bolds{\beta}}{}^{(o)}_{A^*})} {
\|\widehat{\bolds{\beta}}_{A^*}-\widehat{\bolds{\beta}}{}^{(o)}_{A^*}\|} \biggr\rrvert
\\
&&\qquad = \biggl\llvert \biggl(\frac{\partial}{\partial
\bolds{\beta}_{A^*}}Q_n(\widehat{\bolds{\beta}})-\frac{\partial
}{\partial
\bolds{\beta}_{A^*}}Q_n \bigl(\widehat{\bolds{\beta}}{}^{(o)} \bigr) \biggr)^T\frac{(\widehat
{\bolds{\beta}}_{A^*}-\widehat{\bolds{\beta}}{}^{(o)}_{A^*})} {
\|\widehat{\bolds{\beta}}_{A^*}-\widehat{\bolds{\beta}}{}^{(o)}_{A^*}\|} \biggr
\rrvert
\\
&&\qquad = \bigl|n^{-1} \bigl(\widehat{\bolds{\beta}}_{A^*}-\widehat{
\bolds{\beta} }^{(o)}_{A^*} \bigr)^T
\mathbf{X}_{A^*}^T\mathbf{X}_{A^*} \bigl(\widehat{
\bolds{\beta}}_{A^*}-\widehat{\bolds{\beta}}{}^{(o)}_{A^*}
\bigr)/\bigl\|\widehat{\bolds{\beta}}_{A^*}-\widehat{\bolds{\beta}}{}^{(o)}_{A^*}\bigr\|
\\
&&\hspace*{80pt}{}+ \bigl(\bolds{\eta}_{A^*}-\bolds{\xi}^{(o)}_{A^*}
\bigr)^T \bigl(\widehat{\bolds{\beta}}_{A^*} -\widehat{\bolds{\beta}}{}^{(o)}_{A^*} \bigr)/\bigl\|\widehat{\bolds{\beta
}}_{A^*}-\widehat{\bolds{\beta}}{}^{(o)}_{A^*}\bigr\|\bigr|
\\
&&\qquad \geq \phi_{\min} \bigl(n^{-1}\mathbf{X}_{A^*}^T
\mathbf{X}_{A^*} \bigr)\bigl\|\widehat{\bolds{\beta}}_{A^*}-
\widehat{\bolds{\beta}}{}^{(o)}_{A^*}\bigr\|-2\lambda\sqrt
{qu_n^*}
\\
&&\qquad >  \xi_{\min} \bigl(qu_n^* \bigr)2\lambda\sqrt{qu_n^*}\xi_{\min
}^{-1} \bigl(qu_n^*\bigr)-2\lambda\sqrt{qu_n^*}=0,
\end{eqnarray*}
where the second equality follows from the expression of subgradient,
the second last inequality applies the Cauchy--Schwarz inequality, and
the last inequality follows from the relaxed SRC condition in an
$L_0$-neighborhood of the true model. Thus, this contradicts with the
fact that at least one of the subgradients is zero if
$\widehat{\bolds{\beta}}$ is a local minimizer and the theorem is
proved.
\end{pf*}

\begin{pf*}{Proof of Corollary \ref{col_local}}
It follows directly from Theorem \ref{local}.
\end{pf*}


\begin{supplement}
\stitle{Supplement to ``Calibrating nonconvex penalized regression in ultra-high dimension''}
\slink[doi]{10.1214/13-AOS1159SUPP} 
\sdatatype{.pdf}
\sfilename{aos1159\_supp.pdf}
\sdescription{This supplemental material includes the proofs of Lemmas
\ref{lem1} and \ref{le:low}, and some additional numerical results.}
\end{supplement}

%

\printaddresses

\end{document}